\documentclass{amsart}
\usepackage[all]{xypic}
\usepackage{enumerate}

\def\g{\mathfrak g}      
\def\h{\mathfrak h}      
\def\u{\mathfrak u}      
\def\t{\mathfrak t}      
\def\k{\mathfrak k}      
\def\n{\mathfrak n}      
\def\p{\mathfrak p}      
\def\a{\mathfrak a}      
\def\b{\mathfrak b}      
\def\r{\mathfrak r}      
\def\H{\mathcal H}       

\def\w{\mathbf w}        
\def\bu{\mathbf u}       
\def\bl{\mathbf l}       
\def\ba{\mathbf a}       
\def\bL{\mathbf L}       
\def\so{\mathfrak{so}}

\def\pr{\mathrm{pr}}     

\def\R{\mathbb R}        
\def\C{\mathbb C}        
\def\Ad{\mathrm{Ad}}     

\def\Re{\mathrm{Re}}  
\def\Im{\mathrm{Im}}  

\def\Iright{\overset{\rightharpoonup}{I}}
\def\Ileft{\overset{\leftharpoonup}{I}}

\newtheorem{theorem}{Theorem}[section]
\newtheorem{lemma}[theorem]{Lemma}

\newtheorem{proposition}[theorem]{Proposition}

\theoremstyle{definition}
\newtheorem*{defn}{Definition}
\newtheorem*{notation}{Notation}

\theoremstyle{remark}
\newtheorem*{remark}{Remark}

\numberwithin{equation}{section}

\begin{document}
\title{Homogeneous Poisson structures on symmetric spaces}
\author{Arlo Caine}
\address{Max Planck Institut f\"ur Mathematik, P.O. Box: 7280, 53072 Bonn, Germany}
\email{caine@mpim-bonn.mpg.de} 
\author{Doug Pickrell}
\address{Mathematics, University of Arizona, P.O. Box 210089, Tucson, AZ 85721, USA}
\email{pickrell@math.arizona.edu}

\begin{abstract}  We calculate, in a relatively explicit
way, the Hamiltonian systems which arise from the Evens-Lu
construction of homogeneous Poisson structures on both compact and
noncompact type symmetric spaces. A corollary is that the
Hamiltonian system arising in the noncompact case is isomorphic to
the generic Hamiltonian system arising in the compact case.  In
the group case these systems are also isomorphic to those arising
from the Bruhat Poisson structure on the flag space, and hence, by
results of Lu, can be completely factored.
\end{abstract}

\keywords{Lie groups, symmetric spaces, Poisson structures, Hamiltonian systems}


\maketitle

\setcounter{section}{-1}

\section{Introduction}

Suppose that $X$ is a simply connected compact symmetric space
with a fixed basepoint.  From this we obtain  a diagram of groups,
\begin{equation}\label{diagram_of_groups}
\xymatrix{ & G & \\
G_0 \ar[ur] & & U \ar[ul] \\
& K \ar[ul] \ar[ur] & }
\end{equation}
where $U$ is the universal covering of the identity component of
the isometry group of $X$, $X\simeq U/K$, $G$ is the
complexification of $U$, and $X_0=G_0/K$ is the noncompact type
symmetric space dual to $X$; and a diagram of equivariant totally
geodesic (Cartan) embeddings of symmetric spaces:
\begin{equation}\label{CartanEmbeddings}
\xymatrix{
U/K \ar[r]^\phi \ar[d] & U \ar[d] &  \\
G/G_0 \ar[r]^\phi & G & G/U \ar[l]_\psi \\
& G_0 \ar[u] & G_0/K \ar[l]_\psi \ar[u] }
\end{equation}
Let $\Theta$ denote the involution corresponding to the pair
$(U,K)$. We consider one additional ingredient:  a triangular
decomposition of $\g$,
\[
\g=\n^{-}+ \h +\n^{+},
\]
which is $\Theta$-stable and for which $\t_0=\h\cap \k$ is maximal
abelian in $\k$.

This data determines standard Poisson Lie group structures, denoted
$\pi_U$ and $\pi_{G_0}$, for the groups $U$ and $G_0$, respectively.
By a general construction of Evens and Lu (\cite{EL}), the symmetric
spaces $X$ and $X_0$ acquire Poisson structures $\Pi_X$ and
$\Pi_{X_0}$, respectively, which are homogeneous for the respective
actions of the Poisson Lie groups $(U,\pi_U)$ and $(G_0,\pi_{G_0})$.
The compact case was considered in \cite{Caine} and \cite{FL}, and
the noncompact case in \cite{FO}.

In the noncompact case there is just one type of symplectic leaf,
and this leaf is naturally Hamiltonian with respect to the maximal
torus $T_0=\exp(\t_0)\subset K$.  In the compact case the types of
symplectic leaves are indexed by representatives $\w$ of the Weyl
group of $U$ which also lie in the image of the Cartan embedding
$\phi\colon U/K\to U$.  Each such leaf is naturally Hamiltonian
with respect to a torus $T_w$ depending on the corresponding Weyl
group element $w$.  In a reasonably natural way, these leaves are
parameterized by double cosets $R\backslash G_0/K$, where $R$
depends upon $\w$ and a choice of basepoint.  We will refer to these
as Evens-Lu Hamiltonian systems.

The plan of this paper is the following.  In Section
\ref{Background_and_Notation} we recall standard notation  used
throughout the paper.

In Section \ref{Evens-Lu_Hamiltonian_Systems} we exhibit a family
of closed two-forms on $G_0/K$ depending on a parameter $\w_1\in
U$.  For special values of the parameter $\w_1$, these forms
descend to the double coset spaces $R\backslash G_0/K$ and
explicitly describe the Evens-Lu Hamiltonian systems.  The results
in this section for special values of $\w_1$ also follow from the
calculations in sections \ref{The_Noncompact_Case} and
\ref{The_Compact_Case}, and general facts about Poisson geometry.
However this direct approach is suggestive.

In Section \ref{The_Noncompact_Case}, we prove that the
Hamiltonian system with $\w_1=1$ is equal to the Hamiltonian
system arising from the Evens-Lu construction in the noncompact
case.

In Section \ref{The_Compact_Case} we prove that these Hamiltonian
systems are naturally isomorphic to the Hamiltonian systems
arising from the Evens-Lu construction in the compact case. Our
proof of this involves a brutal calculation, which is lacking
conceptual insight.

In Section \ref{The_Group_Case} we specialize to the case
$X=K$, where $K$ is a simply connected compact Lie group.  There
are two main points in this section.  The first is that
$(X,\Pi_X)$ is Poisson isomorphic to the standard Poisson Lie
group structure on $K$, where the isomorphism is essentially
translation by a representative for the longest Weyl group
element.  This translation interchanges the Birkhoff decomposition
(intersected with $K$), the isotypic symplectic components for
$\Pi_X$, with the Bruhat decomposition, the isotypic symplectic
decomposition for the standard Poisson structure.  This
equivalence is a special finite dimensional feature.

The second main point is that the Hamiltonian systems which arise
in this case can all be viewed as torus-invariant symplectic
submanifolds of the generic Hamiltonian system.  This is a
corollary of work of Lu (\cite{Lu}).

In a sequel to this work, we will consider extensions of these
results to loop spaces and applications to the calculation of
integrals.

\section{Background and Notation}\label{Background_and_Notation}

Throughout the remainder of this paper, $U$ will denote a simply
connected compact Lie group, and $\Theta$ will denote an
involution of $U$.  This involution admits a unique holomorphic
extension to $G$ and determines involutive automorphisms of the
Lie algebras of $U$ and $G$, respectively.  Slightly abusing
notation, we will also write $\Theta$ for the extension to $G$ and
the corresponding maps of  algebras.   The identity component of
the fixed point set of $\Theta$ in $U$ will be denoted by $K$, and
$X$ will denote the quotient, $U/K$.  We will also assume that $X$
is irreducible, in the sense of symmetric space theory.

Corresponding to the diagram of groups in
(\ref{diagram_of_groups}), there is a Lie algebra diagram
\begin{equation}\nonumber 
\xymatrix{
& \g = \u + i\u & \\
\g_0 = \k+\p \ar[ur] \ar@/^/[rr]^{\Iright} & & \u = \k+i\p \ar[ul] \ar@/^/[ll]^{\Ileft} \\
& \k  \ar[ul] \ar[ur]& }
\end{equation}
where $\Theta$, acting on the Lie algebra level is $+1$ on $\k$
and $-1$ on $\p$.  The upward arrows are inclusions and the map $\Iright$
(resp. $\Ileft$) is the identity on $\k$ and multiplication by $i$ on $\p$
(resp. $i\p$).  The compositions $\Ileft\circ\Iright$ and $\Iright\circ\Ileft$
agree with $\Theta$ on $\g_0$ and $\u$, respectively.  It will be convenient
to write the action of
$\Theta$ as a superscript, i.e., $\Theta(g)=g^{\Theta}$.  We let
$(\cdot )^{-*}$ denote the Cartan involution for the pair $(G,U)$.
The Cartan involution for the pair $(G,G_0)$ is then given by
$\sigma (g)=g^\sigma=g^{-*\Theta}$.  Since $\Theta$,
$(\cdot)^{-*}$, and $\sigma$ all commute, our practice of writing
these involutions as superscripts should not cause confusion.

There are totally geodesic embeddings of symmetric spaces
\begin{equation}\nonumber
\xymatrix{
U/K \ar[r]^{\phi} \ar[d]& U \ar[d]&\colon & uK \ar[r]& uu^{-\Theta}\\
G/G_0 \ar[r]^{\phi} & G &\colon & gG_0 \ar[r]& gg^{-\sigma} }
\end{equation}
where the symmetric space structures are derived from the Killing
form (the embeddings $\psi$ in (\ref{CartanEmbeddings}) are
defined in a similar way, but will not play a role in this paper).

Fix a maximal abelian subalgebra $\t_0\subset \k$.  By computing
the centralizer $\h_0$ of $\t_0$ in $\g_0$, we obtain
$\Theta$-stable Cartan subalgebras
\begin{equation}\nonumber
\h_0=\t_0+ \a_0,\quad \t=\t_0+ i\a_0,\text{ and } \h=\h_0^{\C}
\end{equation}
for $\g_0$, $\u$, and $\g$, respectively, where $\a_0\subset \p$.
We write
\begin{equation}\label{note_a_0}
\a=\h_{\R}=i\t=i\t_0+\a_0,
\end{equation}
$A=\exp(\a)$, and we let $T_0$ and $T$ denote the maximal tori in
$K$ and $U$ corresponding to $\t_0$ and $\t$, respectively.  We
also fix a $\Theta$-stable triangular decomposition
\[
\g=\n^{-}+ \h +  \n^{+},
\]
so that $\sigma(\n^{\pm})=\n^{\mp}$.  Let
$N^{\pm}=\exp(\n^{\pm})$, $H=\exp(\h)$, $B^{\pm}=HN^{\pm}$, and
let $W$ denote the Weyl group $W(G,H)$. Note that
$W=N_U(T)/T\simeq N_G(H)/H$.

We will write $x=x_{-}+x_{\h}+x_{+}$ for the triangular
decomposition of $x\in \g$, and $x=x_{\k}+x_{\p}$ for the Cartan
decomposition of $x\in \g_0$, and $y=y_\k+y_{i\p}$ for the Cartan
decomposition of $y\in\u$.  To do calculations we will frequently
need to make use of the $\R$-linear orthogonal projections to the
$\R$-subspaces $\u,i\u,\p$, etc.  In keeping with the above
notation scheme, we will write $\{Z\}_\u$ for the orthogonal
projection to $\u$ of $Z\in\g$, and similarly $\{Z\}_{i\u}$,
$\{Z\}_\p$, for the orthogonal projection to $i\u$, $\p$, etc.

There are two decompositions of $\g$ determined by the above data:
\begin{equation}\label{IwasawaDecomps}
\g=\n^-+\a+\u  \text{ and }  \g=\n^-+i\h_0+\g_0.
\end{equation}
The former leads to a global Iwasawa decomposition of the group
$G=N^-AU$. The standard Poisson Lie group structures on $U$ (resp.
$G_0$) that we consider are those associated to the decompositions
in (\ref{IwasawaDecomps}). Given $g\in G$, we write
\begin{equation}\nonumber
g=\bl(g)\ba(g)\bu(g)
\end{equation}
relative to the Iwasawa decomposition $G=N^{-}AU$.   We will write
$\pr_\u$ for the projection $\g\to \u$ with kernel $\n^-+\a$, and
$\pr_{\n^-+\a}$ for the projection $\g\to\n^-+\a$ with kernel
$\u$.  Similarly, $\pr_{\g_0}$ will denote the projection to
$\g_0$ with kernel $\n^-+i\h_0$, and $\pr_{\n^-+i\h_0}$ will
denote the projection to $\n^-+i\h_0$ with kernel $\g_0$.  Note
that the projections $\pr_\u$ (resp. $\pr_{\g_0}$) and the
orthogonal projections $\{\cdot\}_\u$ (resp. $\{\cdot\}_{\g_0}$)
do not agree, as the former has kernel $\n^-+\a$ (resp.
$\n^-+i\h_0$) whereas the latter has kernel $i\u$ (resp. $i\g_0$).

We identify the dual of $\p$ (resp. $i\p$) with $\p$ (resp. $i\p$)
using the Killing form.  To do calculations, we use the induced
isomorphisms
\begin{equation}
G_0\times_K\p\to T(G_0/K)=T^{*}(G_0/K),
\end{equation}
\begin{equation}
U\times_Ki\p\to T(U/K)=T^{*}(U/K).
\end{equation}
To keep track of functoriality, we will write $[g_0,x]$,
$[g_0,y]$, and so on, for tangent vectors, and $[g_0,\phi ]$,
$[g_0,\psi ]$, and so on, for cotangent vectors.

A key player throughout this paper is the ``Hilbert transform''
$\H\colon \g\to\g$ associated to the triangular decomposition of
$\g$,
\begin{equation}\nonumber 
x=x_{-}+x_0+x_{+} \mapsto \H(x)=-ix_{-}+ix_{+}.
\end{equation}
The real subspaces $\g_0$, $i\g_0$, $\u$, and $i\u$ are stabilized
by $\H$, and $\H$ is skew-symmetric with respect to the Killing
form.  This operator also stabilizes $\n^-+\n^+$ and squares to
$-1$ there.  The following proposition is a standard fact about
$\H$.  We include a proof for completeness.

\begin{proposition}\label{vanishingofNtorsionforH}
The Nijenhuis torsion for $\H$ on $\g$,
\begin{equation}\label{NijenhuisforH}
\mathcal N (A,B) = [A,B]+\H([\H(A),B]+[A,\H(B)])-[\H(A),\H(B)],
\quad A,B\in\g,
\end{equation}
 is identically zero.
\end{proposition}
\begin{proof}
Since $\H$ is defined in terms of the triangular decomposition, we
will show that each component of the triangular decomposition of
$\mathcal N (A,B)$ vanishes.  Let $A=A_-+A_0+A_+$ and
$B=B_-+B_0+B_+$ denote the triangular decompositions of $A$ and
$B$.  Since $[A_0,B_0]=0$ we may write
\begin{equation}\label{diagbracket}
[A,B] = [A_-+A_0,B_-+B_0]+([A_-,B_+]+[A_+,B_-])+[A_0+A_+,B_0+B_+].
\end{equation}
The first of the three terms on the right hand side of
(\ref{diagbracket}) is in $\n^-$ since $[\b^-,\b^-]\subset \n^-$.
Similarly, the last term is in $\n^+$.  Hence, the diagonal part
of $[A,B]$ is the same as the diagonal part of the middle term in
(\ref{diagbracket}).  But, $\H$ leaves that term invariant, so we
have $([\H(A),\H(B)])_0=([A,B])_0$.

We can now see that the diagonal of $\mathcal N (A,B)$ vanishes
from the formula in (\ref{NijenhuisforH}).  The second of the
three terms on the right hand side of (\ref{NijenhuisforH}) is in
the image of $\H$ and hence has no diagonal part, whereas we have
just established the equality of the diagonal parts of the first
and last terms.

Let us now turn to the $\n^+$-part of $\mathcal N(A,B)$.  With the
previous observations we have that the $\n^+$ part of the sum of
the first and third terms in the Nijenhuis torsion
(\ref{NijenhuisforH}) is
\begin{equation}\label{Nplusiszero1}
([A,B]-[\H(A),\H(B)])_+=[A_0+A_+,B_0+B_+]+[A_+,B_+].
\end{equation}
Making further use of (\ref{diagbracket}), we compute that
\begin{equation}\label{Nplusiszero2}
([\H(A),B])_+ = (-i[A_-,B_+]+i[A_+,B_-])_+ + i[A_+,B_0+B_+]
\end{equation}
and likewise
\begin{equation}\label{Nplusiszero3}
([A,\H(B)])_+ = (i[A_-,B_+]-i[A_+,B_-])_+ + i[A_0+A_+,B_+].
\end{equation}
Summing the right hand sides of  (\ref{Nplusiszero2}) and
(\ref{Nplusiszero3}) and then applying $\H$ gives that
\begin{equation}\label{Nplusiszero4}
(\H([\H(A),B]+[A,\H(B)]))_+ = - [A_+,B_0+B_+]-[A_0+A_+,B_+]
\end{equation}
which is the $\n^+$-part of the second term in
(\ref{NijenhuisforH}).  The sum of right hand sides of
(\ref{Nplusiszero1}) and (\ref{Nplusiszero4}) gives the
$\n^+$-part of $\mathcal N(A,B)$.
\begin{eqnarray}
(\mathcal N(A,B))_+ & = & [A_0+A_+,B_0+B_+]+[A_+,B_+] \nonumber \\
& & -[A_+,B_0+B_+]-[A_0+A_+,B_+] \label{Nplusiszero5}\\
& = & 0.
\end{eqnarray}
The vanishing of the sum on the right hand side of
(\ref{Nplusiszero5}) is readily apparent after one expands the
terms using bilinearity of the bracket. A completely analogous
calculation shows that the $\n^-$-part of $\mathcal N(A,B)$ is
also zero and completes the proof of Proposition~\ref{vanishingofNtorsionforH}.
\end{proof}

We remark that the equation $\mathcal N(A,B)=0, \forall A,B\in\g$
can be rewritten as
\[
[\H(A),\H(B)]-\H([\H(A),B]+[A,\H(B)])=[A,B]\quad \forall A,B\in\g
\]
which is the modified Yang-Baxter equation for $\H$ with parameter
equal to 1. Using the Killing form, one can view $\H$ as an
element of $\g\wedge\g$, $\u\wedge\u$, or $\g_0\wedge\g_0$.  The
above condition then implies that the Schouten bracket $[\H,\H]$
is $\mathrm{ad}$-invariant as an element of $\wedge^3\g$,
$\wedge^3\u$, or $\wedge^3\g_0$, respectively.  The difference of
the right and left invariant bivector fields generated by $\H$ on
the groups $G_0$ and $U$ determine the standard examples of
Poisson Lie group structures on the semisimple groups $G_0$ and
$U$ (see section 10.4 of \cite{Vaisman}).

Two additional properties of $\H$ which will be important in this
paper concern its relationship with the projections to $\u$ and to
$\g_0$. Given $Z\in\g$, we will write $Z\mapsto iZ$ for the
complex structure on $\g$ and denote the corresponding map of $\g$
by $i$.
\begin{proposition}\label{prop12}
The following diagrams commute.
\begin{equation}\label{Hisprojuionu}
\xymatrix{
\u \ar[r]^{i}\ar[dr]_{\H} & i\u \ar[d]^{\pr_\u} & & \g_0 \ar[r]^i \ar[dr]_{\H} & i\g_0 \ar[d]^{\pr_{\g_0}}\\
            & \u & & & \g_0
}
\end{equation}
\end{proposition}
\begin{proof}
To see that the first diagram commutes, observe that the
triangular decomposition of the element $Z\in \u$ has the form
$-(Z_+)^*+Z_0+Z_+$ where $Z_0\in\t$.  Hence
\[
\pr_\u(iZ)=-(iZ_+)^*+iZ_+=-i(-(Z_+)^*)+iZ_+=\H(Z)
\]
since $i\t=\a$ is contained in the kernel of $\pr_\u$ and the
involution $-(\cdot)^*$ is complex anti-linear.

Similarly, $Z\in \g_0$ has triangular decomposition $(Z_+)^\sigma
+ Z_0+Z_+$ where $Z_0\in \h_0$.  Hence
\[\pr_{\g_0}(iZ) = (iZ_+)^\sigma + iZ_+ = -i(Z_+)^\sigma +iZ_+ = \H(Z)
\] as $i\h_0$ is contained in the kernel of $\pr_{\g_0}$ and the involution $\sigma$ is complex anti-linear.
\end{proof}

\section{Evens-Lu Hamiltonian Systems}\label{Evens-Lu_Hamiltonian_Systems}

In this section we introduce symplectic structures on certain
double coset spaces of $G_0$. The double coset spaces and their
symplectic structures depend on a parameter $\w_1\in U$.  For
certain values of this parameter, these spaces admit Hamiltonian
torus actions for which we compute the momentum maps.

\begin{defn} For each $\w_1\in U$ we define a two-form $\omega_{\w_1}$ on $G_0/K$ by the formula
\begin{equation}\label{defofomega}
\omega_{\w_1} ([g_0,x]\wedge [g_0,y])=\langle
\Ad(\bu(\w_1g_0)^{-1})\circ \H\circ \Ad(\bu(\w_1g_0))(x),y\rangle.
\end{equation}
\end{defn}

\begin{theorem}\label{closureofomegatheorem}
For each $\w_1\in U$, the two-form $\omega_{\w_1}$ on $G_0/K$ is
closed.
\end{theorem}

\begin{proof}
 Fix $\w_1\in U$ and write $\w_1g_0=\bl\ba\bu$ for the Iwasawa factorization of $\w_1g_0\in G$.  Define an $i\u$-valued one-form on $G_0/K$ by $\alpha([g_0,x])=x^\bu$.  Note that the map $g\mapsto \bu(g)$ is right $K$-equivariant, so $\alpha$ is well-defined.  Then $\omega=\langle \H(\alpha)\wedge\alpha\rangle$ and
\begin{equation}\label{211}
d\omega = \langle \H(d\alpha)\wedge\alpha\rangle - \langle \H(\alpha)\wedge d\alpha\rangle.
\end{equation}
Let $\kappa\colon \g_0\to \Gamma(T(G_0/K))$ denote the infinitesimal action of $G_0$ on $G_0/K$.  This is a Lie algebra anti-homomorphism.  Then, given $X,Y\in \g_0$,
\begin{eqnarray}
d\alpha (\kappa(X)\wedge\kappa(Y)) & = & \kappa(X)\alpha(\kappa(Y))-\kappa(Y)\alpha(\kappa(X))-\alpha([\kappa(X),\kappa(Y)]) \nonumber \\
& = & \kappa(X)\alpha(\kappa(Y))-\kappa(Y)\alpha(\kappa(X))+\alpha(\kappa([X,Y])). \label{212}
\end{eqnarray}
To evaluate $d\alpha([g_0,x]\wedge [g_0,y])$ we choose $X$ and $Y$ in $\g_0$ such that the vector fields $\kappa(X)$ and $\kappa(Y)$ agree with the tangent vectors $[g_0,x]$ and $[g_0,y]$ at $g_0K$, respectively.  At $g_0K$, $\kappa(X)$ is represented by $[g_0,\{X^{g_0^{-1}}\}_\p]$, so $X=x^{g_0}$ is one such choice.

The projections $\{\cdot\}_\p$ and $\{\cdot \}_{i\u}$ agree on $\g_0$.  Thus, at $g_0K$,
\[
\alpha(\kappa(X))=(\{X^{g_0^{-1}}\}_{i\u})^\bu=\{X^{\bu g_0^{-1}}\}_{i\u} = \{(X^{\w_1})^{(\bl\ba)^{-1}}\}_{i\u}
\]
since $\w_1 g_0=\bl\ba\bu$ implies that $\bu
g_0^{-1}=(\bl\ba)^{-1}\w_1$.  A straightforward calculation shows
that $(\bl\ba)^{-1}
d(\bl\ba)(\kappa(Y))=\pr_{\n^-+\a}((Y^{\w_1})^{(\bl\ba)^{-1}})$, so
\begin{eqnarray}
\kappa(Y)\alpha(\kappa(X)) & = & \{-\mathrm{ad} (\pr_{\n^-+\a}((Y^{\w_1})^{(\bl\ba)^{-1}}))((X^{\w_1})^{(\bl\ba)^{-1}})\}_{i\u} \nonumber \\
& = & \{-[\pr_{\n^-+\a}(Y^{(\bl\ba)^{-1}\w_1}),X^{(\bl\ba)^{-1}\w_1}]\}_{i\u} \label{213}
\end{eqnarray}
at the point $g_0K$.  Furthermore, we have that
\begin{equation}\label{214}
\alpha(\kappa([X,Y])) =(\{[X,Y]^{g_0^{-1}}\}_{i\u})^\bu =\{[X^{\bu g_0^{-1}},Y^{\bu g_0^{-1}}]\}_{i\u}.
\end{equation}
Thus, with the substitutions $X=x^{g_0}$ and $Y=y^{g_0}$ in
(\ref{213}) and (\ref{214}), we have
$d\alpha([g_0,x]\wedge [g_0,y])=\{W\}_{i\u}$ by (\ref{212}) where
\[
W = -[\pr_{\n^-+\a}(x^\bu),y^\bu]+[\pr_{\n^-+\a}(y^\bu),x^\bu]+[x^\bu,y^\bu].
\]
Now $\pr_\u(\cdot)=\H(-i\cdot)=-i\H(\cdot)$ by Proposition \ref{prop12}.  Thus $\pr_{\n^-+\a}(x^\bu)=x^\bu-\pr_\u(x^\bu)=x^\bu +i\H(x^\bu)$ and
\begin{eqnarray}
W & = & -[x^\bu+i\H(x^\bu),y^\bu]+[y^\bu+i\H(y^\bu),x^\bu]+[x^\bu,y^\bu]\nonumber \\
& = & -i[\H(x^\bu),y^\bu]-i[x^\bu,\H(y^\bu)]-[x^\bu,y^\bu]. \label{215}
\end{eqnarray}
The first two terms of (\ref{215}) are in $i\u$ whereas the third is in $\u$.  Therefore,
\[
d\alpha([g_0,x]\wedge [g_0,y])=\{W\}_{i\u} = -i([\H(x^\bu),y^\bu]+[x^\bu,\H(y^\bu)])
\]
and from (\ref{211}) we have that $d\omega([g_0,x]\wedge [g_0,y]\wedge [g_0,z])$ is equal to
\[
\langle \H (-i([\H(x^\bu),y^\bu]+[x^\bu,\H(y^\bu)])),z^\bu\rangle - \langle \H(x^\bu), -i([\H(y^\bu),z^\bu]+[y^\bu,\H(z^\bu)])\rangle
\]
\[
+\text{ cyclic permutations of }x,y,z.
\]
From the identity $\langle X,[Y,Z]\rangle=\langle [X,Y],Z\rangle$ for each $X,Y,Z\in\g$ and the skew-symmetry of $\H$ with respect to the Killing form, it follows that this sum is equivalent to
\[
i\langle [\H(x^\bu),\H(y^\bu)]-\H([\H(x^\bu),y^\bu]+[x^\bu,\H(y^\bu)]),z^\bu\rangle+i\langle \H(x^\bu),[y^\bu,\H(z^\bu)]\rangle
\]
\[
+\text{ cyclic permutations of }x,y,z.
\]
By (\ref{NijenhuisforH}) the sum of the first, third, and fifth term of the previous expression is equivalent to
\[
i\langle [x^\bu,y^\bu]-\mathcal{N}(x^\bu,y^\bu),z^\bu\rangle + \text{cyclic permutations of }x,y,z
\]
which equals $3i\langle [x^\bu,y^\bu],z^\bu\rangle$ by Proposition \ref{vanishingofNtorsionforH}, whereas the sum of the remaining terms is
\begin{eqnarray}
& & i\langle \H(x^\bu),[y^\bu,\H(z^\bu)]\rangle +i\langle \H(y^\bu),[z^\bu,\H(x^\bu)]\rangle+ i\langle \H(z^\bu),[x^\bu,\H(y^\bu)]\rangle \nonumber \\
& = & i\langle [\H(x^\bu),\H(y^\bu)]-\H([\H(x^\bu),y^\bu]+[x^\bu,\H(y^\bu)]),z^\bu\rangle \nonumber \\
& = & i\langle [x^\bu,y^\bu]-\mathcal{N}(x^\bu,y^\bu),z^\bu\rangle \nonumber \\
& = & i\langle [x^\bu,y^\bu],z^\bu\rangle. \nonumber
\end{eqnarray}
Finally, we have that
\[
d\omega([g_0,x]\wedge [g_0,y]\wedge [g_0,z]) = 4i\langle [x^\bu,y^\bu],z^\bu\rangle = 4i\langle [x,y],z\rangle =0
\]
since $[\p,\p]\subset \k$ which is orthogonal to $\p$. The proof is complete.
\end{proof}

The identification of $U$ with $N^-A\backslash G$ gives rise to a
right action of $G_0$ on $U$.
\[
\begin{matrix}
U\times G_0 & \to & U \\
(u,g_0)& \mapsto & \bu(ug_0)
\end{matrix}
\]
Given $\w_1\in U$, we can compute the stabilizer using the
uniqueness of the Iwasawa decomposition.  Since
$\w_1g_0=g_0^{\w_1}\w_1$, it follows that $\bu(\w_1g_0)=\w_1$ if
and only if $g_0^{\w_1}\in N^-A$.

\begin{notation}
We write $R(\w_1)$ for the $G_0$-subgroup $(N^-A)^{\w_1^{-1}}\cap
G_0$, i.e., the stabilizer of $\w_1$ under the $G_0$-action on
$U$.  The Lie algebra of $R(\w_1)$ will be denoted $\r(\w_1)=(\n^-+\a)^{\w_1^{-1}}\cap\g_0$.
\end{notation}


For the next result concerning $\omega_{\w_1}$ we need a technical device, an operator $\Tilde{\Omega}_0(u)\colon \g_0\to \g_0$ depending on $u\in U$ defined by
\begin{equation}\label{221}
\Tilde{\Omega}_0(u)=\Ileft\circ \Ad(u^{-1})\circ \pr_\u\circ\Ad(u).
\end{equation}
A consequence of part (a) of the following lemma is that if $u=\bu$, the $U$-part of $\w_1g_0=\bl\ba\bu$ as in the proof of the previous theorem, then
\begin{equation}\label{221a}
\omega_{\w_1}([g_0,x]\wedge [g_0,y]) = \langle \Tilde{\Omega}_0(\bu)(x),y\rangle .
\end{equation}

\begin{lemma}\label{lemma22} $\phantom{a}$
\begin{enumerate}
\item[(a)] With respect to the decomposition $\g_0=\k+\p$, $\Tilde{\Omega}_0(u)$ has the form
\begin{equation}\label{222}
\Tilde{\Omega}_0(u)=\begin{pmatrix}1 & b_0 \\ 0 & \Omega_0(u) \end{pmatrix}
\end{equation}
for some linear transformation $b_0\colon \p\to \k$ and where $\Omega_0(u)(x)=\{\Ad(u^{-1})\circ\H\circ\Ad(u)(x)\}_\p$ for each $x\in \p$.
\item[(b)] Suppose that $\w_1\in U$ and $g_0\in G_0$ and write $\w_1g_0=\bl\ba\bu$ for the Iwasawa factorization of $\w_1g_0\in G$.  Then $\Tilde{\Omega}_0(\bu)$ can be factored as the composition
\begin{equation}\label{223}
\xymatrix{
\g_0\ar[r]^{\Ad(g_0)} & \g_0 \ar[r]^{T_{\w_1}(\bu)} &\u \ar[r]^{\Ad(\bu^{-1})} & \u \ar[r]^{\Ileft} &\g_0
}
\end{equation}
where
\begin{equation}\label{224}
T_{\w_1}(\bu)=\pr_\u\circ\Ad((\bl\ba)^{-1}\w_1).
\end{equation}
\item [(c)] Furthermore, $\ker\Omega_0(\bu)=\{\Ad(\r(\w_1))\}_{\p}$.
\end{enumerate}
\end{lemma}

\begin{proof}
Part (a) follows directly from the definition (\ref{221}) together with Proposition \ref{prop12}.  Part (b) is a direct consequence of the factorization $\w_1g_0=\bl\ba\bu$.  For the final claim, note that by (\ref{222}) of part (a), there is an isomorphism $x\in \ker \Omega_0(\bu)\mapsto (-b_0(x),x)\in \ker \Tilde{\Omega}_0(\bu)$ which is section of the orthogonal projection to $\p$.  The factorization in (\ref{223}) implies that $\ker\Tilde{\Omega}_0(\bu)=\Ad(g_0^{-1})(\ker T_{\w_1}(\bu))$.  As an operator on $\g$, the right hand side of (\ref{224}) has kernel
\[
\Ad((\bl\ba)^{-1}\w_1)^{-1}(\n^-+\a)=(\n^-+\a)^{\w_1^{-1}}.
\]
Therefore, $\ker T_{\w_1}(\bu)=(\n^-+\a)^{\w_1^{-1}}\cap \g_0=\r(\w_1)$, and $\ker \Omega_0(\bu)=\{\Ad(g_0)^{-1}(\r(\w_1))\}_{\p}$.  The proof is complete.
\end{proof}


\begin{theorem}
For each $\w_1\in U$, the closed two-form $\omega_{\w_1}$ on $G_0/K$ descends to a
symplectic form on the double coset space $R(\w_1)\backslash
G_0/K$.
\end{theorem}


\begin{proof}
Fix $\w_1\in U$ and for $g_0\in G_0$ write $\w_1g_0=\bl\ba\bu$ for
the Iwasawa factorization of $\w_1g_0\in G$.  It needs to be shown
that the action of $R(\w_1)$ preserves $\omega_{\w_1}$ and $\kappa
(\r(\w_1))|_{g_0K}=\ker(\omega|_{g_0K})$ for each $g_0K\in G_0/K$.
The invariance of $\omega_{\w_1}$ under left translation by
$R(\w_1)$ follows immediately from the definition in
(\ref{defofomega}) given that $R(\w_1)$ is the stabilizer of $\w_1$
in $G_0$ for the right action of $G_0$ on $U$.  From (\ref{221a}) and 
part (a) of Lemma \ref{lemma22},
we know that $\ker(\omega)|_{g_0K}=\ker\Omega_0(\bu)$.  But this is
precisely $\kappa(\r(\w_1))|_{g_0K}$ by part (c) of Lemma
\ref{lemma22}.
\end{proof}


The torus $T_0$ acts on $G_0/K$ by left translation.
In what follows, we introduce other sub-tori of $T$
which will act on the double coset space $R(\w_1)\backslash G_0/K$
for certain values of the parameter $\w_1\in U$.

\begin{notation}
For $w\in W$, we write
\[
\t_w = \{x\in\t\colon \Ad(w)\circ \Theta (x)=x\}
\]
and
\[
T_w = \{t\in T\colon wt^\Theta w^{-1}=t\}.
\]
\end{notation}

Notice that $T_0$ agrees with $T_w$ when $w$ is the trivial element of the Weyl group.

\begin{lemma}\label{momentummaplemma}
Denote by $\w$ the $\w_1K$ Cartan image, $\w=\w_1\w_1^{-\Theta}$.
\begin{enumerate}
\item[(a)] $\Ad(\w_1)\circ\Theta\circ \Ad(\w_1)^{-1}=\Ad(\w)\circ\Theta$
is a complex linear involution of $\g$ which commutes with the
Cartan involution fixing $\u$.
\item[(b)] If $\w\in N_U(T)$ and $w=\w
T\in W$, then:
\begin{enumerate}
\item[(i)] $\h$, $\a$ and $\t$ are $\Theta^{\Ad(\w_1)}$-stable,
\item[(ii)]
$\t_w=\{x\in \t:\Theta^{\Ad(\w_1)}(x)=x\}=\t\cap \g_0^{\w_1}=\t\cap
\k^{\w_1}$, and
\item[(iii)] $T_w=T\cap G_0^{\w_1}=T\cap
K^{\w_1}=\exp(\t_w)$.
\end{enumerate}
\end{enumerate}
\end{lemma}

\begin{proof}
Part (a) follows from the facts that
$\Theta\circ\Ad(\w_1)=\Ad(\w_1^\Theta)\circ\Theta$ and
$\w^{-1}=\w^\Theta$.  Given the validity of (a) and the
$\Theta$-stability of $\h$, $\t$, and $\a$, it follows that each
of these is $\Theta^{\Ad(\w_1)}$-stable when
$\w_1\w_1^{-\Theta}\in N_U(T)$. For (b), part (ii), the set theoretic
description of $\t_w$ follows from (a).  Since $\g_0$ is fixed by
$\sigma$, $\g_0^{\w_1}$ is fixed by $\Ad(\w)\circ \sigma$.  Thus,
by intersecting $\g_0^{\w_1}$ with $\t$ we obtain $\t_w$ which is
the fixed point set in $\t$ of $\Ad(\w)\circ \sigma$.  For the
same reasons, we have that $\t_w=\t\cap\k^{\w_1}$ as $\sigma$ and
$\Theta$ agree and are equal to the identity on $\k$.  The
equalities in (iii) follow routinely from those in (ii).
\end{proof}

\begin{theorem}
Suppose that $\w_1\in U$ is such that $\w=\w_1\w_1^{-\Theta}\in
N_U(T)$ and let $w=\w T$ denote the element of the Weyl group
represented by $\w$.
\begin{enumerate}
\item[(a)] The double coset space $R(\w_1)\backslash G_0/K$ is
contractible.
\item[(b)] The torus $T_w$ acts on $R(\w_1)\backslash
G_0/K$ as follows. Consider $\w_1$ as fixed and abbreviate
$R(\w_1)$ by $R$.
\begin{equation}\label{actionofTw}
\begin{matrix}
T_w \times R\backslash G_0/K & \to & R\backslash G_0/K \\
(t,Rg_0K) & \mapsto & R\w_1^{-1}t\w_1g_0K
\end{matrix}
\end{equation}
Moreover, this action preserves the symplectic form
$\omega_{\w_1}$.
\item[(c)] Let $\t_w^\vee$ denote the dual space of
$\t_w$.  The action of $T_w$ on $R(\w_1)\backslash G_0/K$ is
Hamiltonian with momentum map
\begin{equation}\nonumber
\begin{array}{rccl}
\Phi^{\w_1}\colon & R\backslash G_0/K &  \to & \t_w^\vee \\
 & Rg_0K & \mapsto & \langle i\log \ba(\w_1g_0),\cdot \rangle.
\end{array}
\end{equation}
\end{enumerate}
\end{theorem}

\begin{proof}  Suppose that $\w_1$ and $\w$ are as in the statement of the theorem and regard these as fixed.  We will abbreviate $R(\w_1)$ by $R$ and $\omega_{\w_1}$ by $\omega$ to simplify notation.  For part (a), we refer to the proof of Theorem~4~a) in \cite{pickrell1} which makes use of the assumption $\w\in N_U(T)$.

Let $t\in T_w$.  From Lemma~\ref{momentummaplemma}, we know that
$T_w=T\cap G_0^{\w_1}$. Therefore, $\w_1^{-1}t\w_1\in G_0$ and
$T_w$ acts from the left on $G_0/K$ by $(t,g_0K)\mapsto
t^{\w_1^{-1}}g_0K$.  The fact that
$\Ad(t^{\w_1^{-1}})=\Ad(\w_1^{-1})\circ \Ad(t)\circ\Ad(\w_1)$
preserves $R$ implies that left $T_w$-action on $G_0/K$ descends
to the quotient $R\backslash G_0/K$ as in (\ref{actionofTw}). Note
that
\[
\bu(\w_1(\w_1^{-1}t\w_1)g_0)=t\bu(\w_1g_0)
\]
and $\Ad(t)$ commutes with $\H$.  These observations,
together with the formula in (\ref{defofomega}), imply that
$\omega$ is $T_w$-invariant and proves (b).

Now let $X\in \t_w$, then $\kappa (X^{\w_1^{-1}})$ is a vector field
on $G_0/K$ representing the infinitesimal action of $X$. Let us write
$\w_1g_0=\bl\ba\bu$ for the Iwasawa factorization of $\w_1g_0$.  We must
show that contraction of $\omega$ in the direction of
$\kappa(X^{\w_1^{-1}})$ is equal to the one-form $d\Phi^{\w_1}_X$ where
$\Phi^{\w_1}_X$ is the function
\[
g_0K\mapsto \Phi^{\w_1}_X(g_0K)=\langle iX,\log \ba(\w_1g_0)\rangle .
\]
First we compute $d\Phi^{\w_1}_X$.  Let $\varepsilon$ denote a small real
parameter, let $y\in \p$, and consider $\Phi^{\w_1}_X$ evaluated along
the curve $\varepsilon \mapsto g_0e^{\varepsilon Y}K$.  Observe
that
\[
\left.\frac{d}{d\varepsilon}\right|_{\varepsilon = 0} \log
\ba(\w_1 g_0e^{\varepsilon y}) = [g_0,\{y^\bu\}_{\a}]
\]
since the orthogonal projection to $\a$ and the Iwasawa projection
to $\a$ give the same result.  Thus
\[
d\Phi^{\w_1}_X ([g_0,y]) = \langle iX, \{y^\bu\}_\a\rangle = \langle iX,
y^\bu\rangle = \langle iX^{\bu^{-1}},y\rangle
\]
since $X\in \t_w\subset \t$.  This shows that $d\Phi^{\w_1}_X$ is
represented by the class $[g_0, \{X^{\bu^{-1}}\}_\p]$. 

At $g_0K$, $\kappa(X^{\w_1^{-1}})=[g_0,\{X^{(\w_1 g_0)^{-1}}\}_\p]$. Part (b) (ii) of Lemma~\ref{momentummaplemma} implies that $\Ad(\w_1g_0)^{-1}(X)\in\g_0$
and therefore $\{\Ad(\w_1g_0)^{-1}(X)\}_\p=\{\Ad(\w_1g_0)^{-1}(X)\}_{i\u}$.
Hence, we have that
\begin{eqnarray}
\omega(\kappa(X^{\w_1^{-1}})\wedge [g_0,y]) & = & \langle \Ad(\bu^{-1})\circ \H\circ \Ad(\bu)(\{\Ad((\w_1g_0)^{-1})(X)\}_\p),y\rangle \nonumber \\
& = & \langle \H\circ \Ad(\bu)(\{\Ad((\w_1g_0)^{-1})(X)\}_{i\u}), y^\bu\rangle  \nonumber \\
& = & \langle \{\H(X^{(\bl\ba)^{-1}})\}_{i\u}, y^\bu \rangle \nonumber \\
& = & \langle \{-i(X^{(\bl\ba)^{-1}})_-\}_{i\u}, y^\bu \rangle \label{Thm24_A}
\end{eqnarray}
using the factorization $\w_1g_0=\bl\ba\bu$ and the fact that $\H$ commutes with the orthogonal
projection to $i\u$.  To continue, note that $(X^{(\bl\ba)^{-1}})_-=X^{(\bl\ba)^{-1}}-X$ because
$X\in \t_w\subset \t$ and thus
$X^{(\bl\ba)^{-1}}\in \b^-$.  Continuing from (\ref{Thm24_A}), we have that
\begin{eqnarray}
\omega (\kappa (X^{\w_1^{-1}})\wedge [g_0,y]) & = & \langle \{-i(X^{(\bl\ba)^{-1}})_-\}_{i\u}, y^\bu \rangle \nonumber \\
& = & \Im\langle \{(X^{(\bl\ba)^{-1}})_-\}_{\u}, y^\bu \rangle \nonumber \\
& = & \Im\langle (X^{(\bl\ba)^{-1}})_-, y^\bu \rangle \nonumber \\
& = & \Im\langle X^{(\bl\ba)^{-1}}, y^\bu \rangle - \Im\langle X,y^\bu \rangle \label{Thm24_B}
\end{eqnarray}
using that multiplication by $i$ intertwines the orthogonal projections to $i\u$ and $\u$ and that
the Killing form is complex linear and real valued on $i\u$.  The first term on the right hand side
of (\ref{Thm24_B}) vanishes because
\[
\Im\langle X^{(\bl\ba)^{-1}}, y^\bu \rangle = \Im\langle X, y^{\w_1g_0}\rangle
\]
and both factors are in $\g_0^{\w_1}$ which is a real form of $\g$.  The second term on
the right hand side of (\ref{Thm24_B}) can be rewritten $ - \Im\langle X,y^\bu \rangle=\langle iX,y^\bu\rangle$.
Thus, we have shown that
\[
\omega (\kappa (X^{\w_1^{-1}})\wedge  [g_0,y]) = \langle
iX,y^\bu\rangle = \langle iX^{\bu^{-1}},y\rangle = d\Phi^{\w_1}_X
([g_0,y])
\]
completing the proof of the theorem.
\end{proof}

\section{The Noncompact Case}\label{The_Noncompact_Case}

We will write $X_0$ for the non-compact symmetric space $G_0/K$.
The Evens-Lu Poisson structure on $X_0$ is given by the formula
\begin{equation}\nonumber 
\Pi_{X_0}([g_0,\phi ]\wedge [g_0,\psi
])=\langle\Omega(g_0)(\phi),\psi\rangle ,
\end{equation}
where
\begin{equation}\nonumber
\Omega (g_0)(\phi )=\{\Ad(g_0)^{-1}\circ \H\circ \Ad(g_0)(\phi)\}_{\p}.
\end{equation}
Note that $\Omega(g_0)\in\so(\p)$ because $\H$ is skew.

 In this section, we
explicitly describe the geometry of the symplectic foliation for
$\Pi_{X_0}$.  This structure is regular and we can compute a
Casimir.  Lastly, we show that along the symplectic leaves in
$G_0/K$, the two-form $\Pi_{X_0}^{-1}$ agrees with the restriction
of the global two-form $\omega_{\w_1}$ (introduced in section
\ref{Evens-Lu_Hamiltonian_Systems}) with the parameter $\w_1=1\in
U$.

To obtain the results we need to work with an extension of
$\Omega(g_0)$ to all of $\g_0$ defined by
\begin{equation}\label{dougsextendedoperator}
\Tilde{\Omega}(g_0)=\Ad(g_0)^{-1}\circ \pr_{\g_0}\circ \Ad(g_0)\circ
\Iright.
\end{equation}
A consequence of part (a) of the following lemma is that
\begin{equation}\nonumber
\Pi_{X_0}([g_0,\phi]\wedge [g_0,\psi])=\langle
\Tilde{\Omega}(g_0)(\phi),\psi\rangle.
\end{equation}

\begin{lemma}\label{lemma30}$\phantom{a}$
\begin{enumerate}
\item[(a)] Relative to
the decomposition $\k+\p$,
\begin{equation}\label{presentationofOmega}
\Tilde{\Omega}(g_0)=\begin{pmatrix}
1 & b \\
0 & \Omega(g_0)
\end{pmatrix}
\end{equation}
for some linear transformation $b\colon \p\to\k$ depending on $g_0$.
\item[(b)] Write $g_0=\bl\ba\bu=\bl\ba_0\ba_1\bu$ for the Iwasawa
factorization of $g_0$ with the further factorization
$\ba=\ba_0\ba_1$ where $\ba_0$ is the $A_0$ part of $\ba$ (cf.
(\ref{note_a_0})).  Then $\Tilde{\Omega}(g_0)$ can be factored as
the composition
\begin{equation}\label{factorizationofOmega}
\xymatrix{ \g_0 \ar[r]^{\Iright} & \u \ar[r]^{\Ad(\bu)} & \u
\ar[r]^{T(g_0)} & \g_0 \ar[rr]^{\Ad((\ba_0^{-1}g_0))^{-1}}&  & \g_0,
}
\end{equation}
where for $X\in \u$,
\begin{equation}\label{formula_for_T}
T(g_0)(X)=(((X_+)^\bL)_+)^\sigma +
X_{\t_0}+((X_+)^\bL)_{\h_0}+((X_+)^\bL)_+
\end{equation}
with $\bL=\ba_0^{-1}\bl\ba_0\ba_1$.
\end{enumerate}
\end{lemma}

\begin{proof}
Part (a) follows from the definition and Proposition \ref{prop12}.
For (\ref{factorizationofOmega}), set $T(g_0)=\Ad(\ba_0^{-1})\circ
\pr_{\g_0}\circ \Ad(\bl\ba)$.  Because $\ba_0\in G_0\cap H$,
$\Ad(\ba_0)$ commutes with $\sigma$ and stabilizes the triangular
decomposition.  Hence it also commutes with $\pr_{\g_0}$, so
$T(g_0)=\pr_{\g_0}\circ \Ad(\bL)$ where
$\bL=\ba_0^{-1}\bl\ba_0\ba_1$.  The formula for $T(g_0)(X)\in \g_0$
then follows from the identities $(X^\bL)_+=((X_+)^\bL)_+$ and
$(X^\bL)_{\h_0}=X_{\h_0}+((X_+)^\bL)_{\h_0}=X_{\t_0}+((X_+)^\bL)_{\h_0}$
which we have because $X\in\u$ and $\L\in B^-$.
\end{proof}

\begin{remark}
The formula
$\Omega(g_0)(\phi)=\{(\pr_{\g_0}(i\phi^{g_0}))^{g_0^{-1}}\}_{\p}$ we
get from part (a) of Lemma \ref{lemma30}, first appeared in
\cite{pickrell1} along with the extended operator
$\Tilde{\Omega}(g_0)$. See section 5 of that paper for a derivation
of this formula from the general construction in \cite{EL}. In
\cite{pickrell1}, displayed line (37), $\bl'$ should have been set
equal to $\bl'=\ba_0^{-1}\bl\ba_0$ rather than $\ba_0\bl\ba_0^{-1}$.
However, this has no effect on the remaining results in that paper.
\end{remark}

For later purposes, we now establish a number of facts about the
operator $T(g_0)$.  Write $\u\ominus i\a_0$ for the orthogonal complement of $i\a_0$ in $\u$ with respect to the Killing form.

\begin{lemma}\label{PropofT}For each $g_0\in G_0$ write $g_0=\bl\ba_0\ba_1\bu$ and set $\bL=\ba_0^{-1}\bl\ba_0\ba_1$ as in (b) of Lemma \ref{lemma30}.
\begin{enumerate}
\item[(a)] $\ker(T(g_0))=i\a_0$.
\item[(b)] The adjoint of $T(g_0)$, relative to the Killing forms on $\u$ and $\g_0$, $T^*(g_0)\colon \g_0\to\u$ is given by the formula
\begin{equation}\label{formula_for_Tadjoint}
T^*(g_0)(y)=\frac{1}{2}\left(((y_{\h_0}+2y_-)^{\bL^{-1}})_-+2y_{\t_0}-(((y_{\h_0}+2y_-)^{\bL^{-1}})_-)^*\right).
\end{equation}
\item[(c)] The cokernel of $T(g_0)$ is
$\ker(T^*(g_0))=\left\{(y^\bL-y)+2y+(y^\bL-y)^\sigma\colon y\in
\a_0\right\}.$
\item[(d)] The image of $T(g_0)$ consists of $y\in \g_0$ such that the
$\h_0$ part of $y$ is in the image of the following map.
\begin{equation}\nonumber
\begin{array}{rcl} \t_0+\n^+ & \to & \h_0 \\
x_{\t_0}+x_+ & \mapsto & x_{\t_0}+\{(x_+)^\bL\}_{\h_0}
\end{array}
\end{equation}
Given such a $y\in\g_0$, the solution of $T(x)=y$, for $x\in \u\ominus i\a_0$, is solved in stages by
\begin{equation}\nonumber
x_+=((y_+)^{\bL^{-1}})_+\text{, }x_- =-(x_+)^*, \text{ and }
x_{\t_0}=y_{\t_0}-\{(x_+)^\bL\}_{\t_0}.
\end{equation}
The moral is that it is easy to solve for $T(g_0)^{-1}$ if one
knows that there is a solution.
\end{enumerate}
\end{lemma}

\begin{proof}
As an operator on $\g$, $\pr_{\g_0}\circ \Ad(\bL)$ has kernel
$\Ad(\bL^{-1})(\n^-+i\h_0)=\n^-+i\h_0$ because $\bL\in B^-$. This
establishes (a) because $\ker(T(g_0))=\ker(\pr_{\g_0}\circ
\Ad(\bL))\cap\u=(\n^-+i\h_0)\cap\u = i\a_0$.

 For part (b), we first show that $\Re\langle
\pr_{\g_0}(Z),y\rangle =\langle Z, 2y_-+y_{\h_0}\rangle$ for each
$Z\in \g$ and $y\in \g_0$.  Indeed,
\begin{eqnarray}
\Re\langle \pr_{\g_0}(Z),y\rangle & = &
\Re\langle(Z_+)^\sigma, (y_-)^\sigma\rangle +\Re\langle\{Z\}_{\h_0},y_{\h_0}\rangle +\Re\langle Z_+,y_-\rangle \label{311}\\
& = & \Re\langle Z_+,y_-\rangle+\Re\langle Z_0,y_{\h_0} \rangle
+\Re\langle Z_+,y_-\rangle \label{312}\\
& = & \Re\langle Z,2y_-+y_{\h_0}\rangle.\nonumber
\end{eqnarray}
Note that (\ref{312}) follows from (\ref{311}) only because we were using
the real part of the Killing form.  This claim with $Z=\Ad(\bL)(x)$
for $x\in \u$, gives
\begin{eqnarray*}
\langle T(g_0)(x),y \rangle & = & \Re\langle \pr_{\g_0}\circ
\Ad(\bL)(x),y\rangle \\
& = & \Re\langle x, \Ad(\bL^{-1})(2y_-+y_{\h_0})\rangle,
\end{eqnarray*}
and hence $T^*(g_0)(y)=\{\Ad(\bL^{-1})(2y_-+y_{\h_0})\}_\u$ by
non-degeneracy of the Killing form (which is real valued) on $\u$.
The zero mode of $T^*(g_0)(y)$ is then
$\{(\Ad(\bL)^{-1}(2y_-+y_{\h_0}))_0\}_\u=\{y_{\h_0}\}_\u=y_{\t_0}$.
The formula in (\ref{formula_for_Tadjoint}) follows immediately.

Now suppose that $T^*(g_0)(y)=0$.  By (\ref{formula_for_Tadjoint}),
$y_{\t_0}=0$, so $y_{\h_0}=y_{\a_0}$, and
$((y_{\a_0}+2y_-)^{\bL^{-1}})_-=0$.  We can use this last equation
to determine $y_-$ in terms of $y_{a_0}$.  After all,
\[
0=((y_{\a_0}+2y_-)^{\bL^{-1}})_-
=((y_{\a_0})^{\bL^{-1}})_-+((2y_-)^{\bL^{-1}})_- =
(y_{\a_0})^{\bL^{-1}}-y_{\a_0}+(2y_-)^{\bL^{-1}},
\]
whence $2y_- = y_{\a_0}-y_{\a_0}^{\bL}$.  After rescaling $y$ to
$2y$ we obtain the description of the elements of the cokernel in
(c).

The first part of (d) concerning the image of $T(g_0)$ follows
easily after examining the formula (\ref{formula_for_T}).  Let
$x\in\u\ominus i\a_0$ and $y$ be in the image of $T(g_0)$.  Then
\begin{equation}\nonumber
y_+ = ((x_+)^\bL)_+ = (x_+)^\bL + z
\end{equation}
where $z\in\b^-$.  Using that $\bL^{-1}\in B^-$ we obtain that
$((y_+)^{\bL^{-1}})_+=x_+$. Once $x_+$ is determined, we know that
we can find $x_{\t_0}$ by the equation
\begin{equation}\nonumber
x_{\t_0}=y_{\h_0}-\{(x_+)^\bL\}_{\h_0}=y_{\t_0}-\{(x_+)^\bL\}_{\t_0}.
\end{equation}
This completes the proof of the Lemma~\ref{PropofT}.
\end{proof}

\begin{lemma}\label{tangenttoleafcondition}$\phantom{a}$
\begin{enumerate}
\item[(a)] The tangent vector $[g_0,x]$ is tangent to the symplectic leaf
through $g_0K$ if and only if $\Ad(\bu)(x)$ is perpendicular to
$\a_0$ relative to the Killing form on $i\u$.
\item[(b)] The Poisson structure $\Pi_{X_0}$ is regular.
\end{enumerate}
\end{lemma}

\begin{proof}
As usual, write $g_0=\bl\ba\bu$ for the Iwasawa factorization of
$g_0$.  The subspace tangent to the symplectic leaf through $g_0K$
is the image of the anchor map $\Pi_{X_0}^\#\colon T^*(G_0/K)\to
T(G_0/K)$ at $g_0K$.  In terms of our working identifications
\[
\Pi_{X_0}^\#([g_0,\phi])=[g_0,\Omega(g_0)(\phi)].
\]
Since $\Omega(g_0)\in \mathfrak{so}(\p)$ for each $g_0\in G_0$, its
image is equal to the orthogonal complement of its kernel. It
follows from (\ref{presentationofOmega}) that there is an
isomorphism $\phi\in \ker(\Omega(g_0)) \mapsto (-b(\phi),\phi)\in \ker
(\Tilde{\Omega}(g_0))$ which is a section of the projection to $\p$.
The factorization of $\Tilde{\Omega}(g_0)$ in
(\ref{factorizationofOmega}) together with part (a) of
Lemma~\ref{PropofT} shows that
\[
\ker(\Tilde{\Omega}(g_0))=\Iright\hspace{-0.5em}\phantom{a}^{-1}\circ\Ad(\bu)^{-1}(\ker
(T(g_0)))=\Iright\hspace{-0.5em}\phantom{a}^{-1}\circ\Ad(\bu^{-1})(i\a_0).
\]
Thus $\ker \Omega(g_0)=
-i\{\Ad(\bu^{-1})(i\a_0)\}_{i\p}=\{\Ad(\bu^{-1})(\a_0)\}_\p$ and the
map $\a_0\to \{\Ad(u^{-1})(\a_0)\}_\p$ is an isomorphism for each
$u\in U$ such that $u=\bu(g_0)$ for some $g_0\in G_0$. This proves
Lemma \ref{tangenttoleafcondition}.
\end{proof}

\begin{proposition}$\phantom{a}$
\begin{enumerate}
\item[(a)] The
image of the anchor map $\Pi_{X_0}^\#\colon T^*(G_0/K)\to
T(G_0/K)$ defines a flat connection for the principal bundle
\[
\xymatrix{A_0 \ar[r] & G_0/K \ar[r] & A_0\backslash G_0/K.}
\]
\item[(b)] The symplectic leaves are the level sets of the function
$\ba_0$.
\item[(c)] The horizontal parameterization for the symplectic
leaf through the basepoint is given by the map $s\colon
A_0\backslash G_0/K\to G_0/K$
\begin{equation}\label{leafsection}
A_0 g_0 K \to s(A_0g_0K)=\ba_0^{-1}g_0K
\end{equation}
where $g_0=\bl\ba_0\ba_1\bu$.
\end{enumerate}
\end{proposition}

\begin{proof} We continue to write $g_0=\bl\ba\bu$
for the Iwasawa factorization of $g_0$. Essentially, part (a) (along
with part (b) of the previous lemma) was established in \cite{FO}.
We supply an alternative argument here.

Given (b) of Lemma \ref{tangenttoleafcondition}, it suffices to
check infinitesimally that the $A_0$-orbits have trivial
intersection with the symplectic leaves. The connection is flat
because the symplectic leaf distribution is integrable. The tangent
space to the $A_0$-orbit through $g_0K$ is
$\{[g_0,\{\Ad(g_0^{-1})(y_0)\}_\p]\colon y_0\in \a_0\}$. It is clear
that at the basepoint this subspace intersects the subspace tangent
to the symplectic leaf only at zero.  In general, let $y_0\in\a_0$
and suppose that there exists $x\in\p$ with $x^\bu \perp \a_0$ such
that $\{\Ad(g_0^{-1})(y_0)\}_\p=x$.  Let $\kappa\in \k$ be such that
$\Ad(g_0^{-1})(y_0)=\kappa + x$. Given $z_0\in \a_0$, the pairing
$\langle \Ad(\bu g_0^{-1})(y_0),z_0\rangle$ can be written in two
equivalent ways. On the one hand $\Ad(\bu
g_0^{-1})(y_0)=\Ad((\bl\ba)^{-1})(y_0)$ so
\begin{equation}\label{A0orbits1}
\langle \Ad(\bu g_0^{-1})(y_0),z_0\rangle = \langle
\Ad((\bl\ba)^{-1})(y_0) , z_0\rangle = \langle y_0 , z_0\rangle
\end{equation}
since $y_0\in \h$ and $\bl\ba\in B^-$.  On the other hand
\begin{equation}\label{A0orbits2}
\langle \Ad(\bu g_0^{-1})(y_0),z_0\rangle = \langle \kappa^\bu +
x^\bu , z_0\rangle =\langle \kappa^\bu , z_0\rangle
\end{equation}
since $x^\bu \perp \a_0$.  The right hand side of
(\ref{A0orbits1}) is real whereas the right hand side of
(\ref{A0orbits2}) is purely imaginary since $\kappa^\bu\in \u$, so
they must both be zero.  This implies that $x$ must be zero,
proving (a).

For part (b), identify the tangent bundle to $A_0$ with $A_0\times
\a_0$ using left translation. Then $d\ba_0$ is identified with the
$\a_0$-valued one-form $[g_0,x]\mapsto \{x^\bu\}_{\a_0}$.  It then
follows from Lemma~\ref{tangenttoleafcondition} that the
symplectic leaves are the level sets of $\ba_0$.

We now turn to part (c).  Let $a_0\in A_0$, then $a_0g_0=a_0\bl\ba_0\ba_1\bu = \bl^{a_0}a_0\ba_0\ba_1\bu$.
Since $\Ad(a_0)$ stabilizes $N^-$, it follows from the uniqueness
of the Iwasawa decomposition that the $A_0$ factor of $a_0g_0$ is
$a_0\ba_0$.  This shows that the cross section (\ref{leafsection})
is well-defined.  It remains to show that the image of $s$ is
horizontal.

Let $\varepsilon$ be a small real parameter so that, given
$x\in\p$, the map $\varepsilon \mapsto g_0e^{\varepsilon x}K$ is a
smooth curve passing through $g_0K$ at $\varepsilon = 0$.  Then
\begin{equation}\nonumber
\left.\frac{d}{d\varepsilon}\right|_{\varepsilon = 0}
s(g_0e^{\varepsilon x}K) =[\ba_0^{-1}g_0, x -
\{(\{x^\bu\}_{\a_0})^{g_0^{-1}}\}_\p].
\end{equation}
To show that this is horizontal, we must check that the pairing
\begin{equation}\label{sectionishorizontal1}
\langle (x - \{(\{x^\bu\}_{\a_0})^{g_0^{-1}}\}_\p)^\bu, y_0
\rangle
\end{equation}
vanishes for $y_0\in \a_0$.  Again, the projection to $\p$ in
(\ref{sectionishorizontal1}) may be replaced with the projection
to $i\u$ as it is being applied to an element of $\g_0$.  Thus,
(\ref{sectionishorizontal1}) becomes
\begin{equation}\label{sectionishorizontal2}
\langle x^\bu, y_0\rangle  - \langle
\{(\{x^\bu\}_{\a_0})^{(\bl\ba)^{-1}}\}_{i\u}, y_0 \rangle
\end{equation}
using the factorization $g_0=\bl\ba\bu$ and the fact that $\Ad(\bu)$
commutes with the projection to $i\u$.  Note that if $Z\in
i\t\subset i\u$, then the diagonal part of
$\{\Ad((\bl\ba)^{-1})(Z)\}_{i\u}$ is $Z$ since $\bl\ba\in B^-$.
Apply this observation to $Z=\{x^\bu\}_{\a_0}$ and we have that the
second term in (\ref{sectionishorizontal2}) is
\[
\langle \{(\{x^\bu\}_{\a_0})^{(\bl\ba)^{-1}}\}_{i\u}, y_0 \rangle
= \langle \{x^\bu\}_{\a_0},y_0\rangle = \langle x^\bu ,
y_0\rangle.
\]
Therefore (\ref{sectionishorizontal1}) vanishes.  The proof is
complete.
\end{proof}

\begin{theorem}
Along the symplectic leaves, $\Pi_{X_0}^{-1}$ agrees with the
restriction of the closed two-form $\omega_{\w_1}$ from
(\ref{defofomega}) with $\w_1=1$.
\end{theorem}

\begin{proof}
Factor $g_0\in G_0$ as $\bl\ba_0\ba_1\bu$ as before and set
$\bL=\ba_0^{-1}\bl\ba_0\ba_1\in B^-$, then $\ba_0^{-1}g_0=\bL\bu$.
Let $[g_0,x]$ and $[g_0,y]$ represent tangent vectors to the
symplectic leaf through $g_0K$.  Making use of the extended operator
(\ref{dougsextendedoperator}) we have
\begin{equation}\label{341}
\Pi_{X_0}^{-1}([g_0,x]\wedge [g_0,y])=\langle
\Omega^{-1}(g_0)(x),y\rangle = \langle
\Tilde{\Omega}^{-1}(g_0)(x),y\rangle
\end{equation}
where the inverses of $\Omega(g_0)$, $\Tilde{\Omega}(g_0)$ and
$T^{-1}(g_0)$ are computed on the orthogonal complement on their
kernels. The factorization (\ref{factorizationofOmega}), with the
substitution $\ba_0^{-1}g_0=\bL\bu$, implies that
\[
\Tilde{\Omega}^{-1}(g_0) = \Iright\hspace{-0.5em}\phantom{a}^{-1}\circ \Ad(\bu)^{-1}\circ
T^{-1}(g_0)\circ \Ad(\bL)\circ\Ad(\bu).
\]
Our goal is compute $\Tilde{\Omega}^{-1}(g_0)(x)$ and to that end we
will first  compute $X\in \u$ such that
$T(g_0)(X)=\Ad(\bL)\circ\Ad(\bu)(x)$.  Set $\chi=\Ad(\bu)(x)=x^\bu$.
By part (a) of Lemma \ref{tangenttoleafcondition}, $\chi$ is orthogonal to
$\a_0$, so
\[
(\Ad(\bL)(x))_{\h_0}=(\chi_0)_{\h_0}+(((\chi_+)^\bL)_0)_{\h_0}=((\chi_+)^\bL)_{\h_0}
\]
since $\chi_0\in i\t_0$.  It now follows from Lemma \ref{PropofT},
part (d), that  there exists $X\in\u$ such that
$T(g_0)(X)=\Ad(\bL)(\chi)$.  Furthermore,
$X_+=(((\chi^\bL)_+)^{\bL^{-1}})+=\chi_+$ since $L\in B^-$,
$X_-=-(X_+)^*$, $X_{i\a_0}=0$, and
\[
X_{\t_0}=\{\chi^\bL\}_{\t_0} -\{(\chi_+)^\bL\}_{\t_0}=\{\chi_0\}_{\t_0}+\{(\chi_+)^\bL\}_{\t_0} -\{(\chi_+)^\bL\}_{\t_0}=0.
\]
Thus, $X=-(\chi_+)^*+\chi_+=-(x^\bu)_-+(x^\bu)_+$ because
$\chi=x^\bu\in i\u$.  We now have from (\ref{341})
\begin{eqnarray*}
\Pi_{X_0}^{-1}([g_0,x]\wedge [g_0,y]) & = & \langle \Iright\hspace{-0.5em}\phantom{a}^{-1}\circ \Ad(\bu)^{-1}(-(x^\bu)+(x^\bu)_+),y\rangle \\
& = & -\langle \Ad(\bu)\circ\H\circ\Ad(\bu)(ix),iy\rangle \\
& = & \omega_{\w_1}([g_0,x]\wedge [g_0,y])
\end{eqnarray*}
with $\w_1=1\in U$.
\end{proof}

\section{The Compact Case}\label{The_Compact_Case}

The Evens-Lu Poisson structure on $X=U/K$ is given by the formula
\begin{equation}\nonumber
\Pi_X ([u,\phi]\wedge[u,\psi])=\langle \Omega(u)(\phi),\psi\rangle
\end{equation}
where the linear transformation $\Omega(u)\colon i\p\to i\p$ is
given by
\begin{equation}\nonumber
\Omega (u)(\phi) = \{\Ad(u)^{-1}\circ \H\circ\Ad(u)(\phi)\}_{i\p}.
\end{equation}
Note that $\Omega(u)\in\so(i\p)$ because $\H$ is skew.  See section 2 of \cite{Caine} for a derivation of this formula
from the Evens-Lu construction. Recall that $G_0$ acts from the
right on $U$ through the Iwasawa decomposition.
\[
\begin{matrix}
U\times G_0 & \to & U\\
(u,g_0) & \mapsto & \bu(ug_0)
\end{matrix}
\]
It was shown in \cite{FL} that the symplectic leaves of $\Pi_X$
are the projections of the $G_0$-orbits in $U$ to $U/K$.  Building
on this work and that of \cite{pickrell1}, a finer description was
given in \cite{Caine} using the connection with the Birkhoff
decomposition.

Corresponding to the triangular decomposition $\g = \n^- +\h+
\n^+$ there is a decomposition of the group
\[
G=\coprod_{w\in W} \Sigma_w, \text{ where the }\Sigma_w = N^-wHN^+
\]
are submanifolds whose complex codimension is equal to the
length of the indexing Weyl group element.  The symmetric space
$X$ inherits a decomposition into the pre-images of the $\Sigma_w$
under the map
\[
\xymatrix{ X=U/K \ar[r] & U \ar[r] & G }
\]
where the first arrow is the Cartan embedding and the second is
inclusion.  As a variety in $U$, the image of the Cartan embedding
is the connected component containing the identity of
$\{u^{-1}=u^\Theta\}\subset U$.  As in \cite{Caine}, the pre-image
of $\Sigma_w$ will be referred to as the layer of the Birkhoff
decomposition indexed by $w$.  Literally viewing the Weyl group
$W=N_U(T)/T$ as the set of connected components of the normalizer
of $T$ in $U$, one obtains that the layers of the Birkhoff
decomposition of $X$ are indexed by those elements $w\in W$ such
that $w\cap \{u^{-1}=u^\Theta \}_0\not=\emptyset$.  Each layer may
consist of multiple connected components.

Because the triangular decomposition of $\g$ is $\Theta$-stable,
the symplectic foliation of $\Pi_X$ aligns with the Birkhoff
decomposition.  Each connected component of a given layer is
foliated by contractible symplectic leaves. When restricted to a
given layer, $\Pi_X$ is regular. The torus $T_w=\{t\in T\colon
wt^\Theta w^{-1} =t\}$ (cf. section
\ref{Evens-Lu_Hamiltonian_Systems}) acts on the layer indexed by
$w$ preserving the symplectic leaves.  The action on each leaf is
Hamiltonian and has a unique fixed point.  The images of the
$T_w$-fixed points under the Cartan embedding are the elements of
the intersection of the image of the Cartan embedding with
$w\subset U$. We thus label the symplectic leaves of $(X,\Pi_X)$
by the representatives $\w\in w \cap \{u^{-1}=u^\Theta \}_0$.
\begin{notation}
 We will denote by $S(\w)$ the symplectic leaf of $(X,\Pi_X)$ corresponding to $\w$. When we write, ``Let $S(\w)$ be a symplectic leaf,'' we implicitly declare that $\w$ is in $N_U(T)$ and in the image of the Cartan embedding.  By $\Pi_\w$ we denote the restriction of the Poisson tensor $\Pi_X$ to the symplectic leaf $S(\w)$.
\end{notation}

Let $S(\w)$ be a symplectic leaf.  Fix a choice of $\w_1\in U$
such that $\w_1\w_1^{-\Theta}=\w$. The map $\Tilde{\bu}\colon G_0
\to U$ defined by $g_0\mapsto \Tilde{\bu}(g_0)=\bu(\w_1g_0)$
is equivariant for the right actions of $K$ on $G_0$ and $U$,
invariant under the left action of $R(\w_1)=(N^-A)^{\w_1^{-1}}\cap
G_0$ on $G_0$, and descends to a $T_w$-equivariant diffeomorphism
\begin{equation}\nonumber 
\Tilde{\bu}\colon  R(\w_1)\backslash G_0/ K  \to  S(\w).
\end{equation}

The main result of this section is the following theorem.

\begin{theorem}\label{isomtheorem} Let $S(\w)$ be a symplectic leaf.  Fix a choice of $\w_1\in U$
such that $\w_1\w_1^{-\Theta}=\w\in N_U(T)$.  Then the map
$\Tilde{\bu}$ induces an isomorphism of $T_w$-Hamiltonian spaces
\begin{equation}
(R(\w_1)\backslash G_0/K,\omega_{\w_1})\to (S(\w),\Pi_{\w}^{-1})
\end{equation}
where $\omega_{\w_1}$ is as in (\ref{defofomega}).
\end{theorem}

We remark here that there is a sense in which this result does not
depend upon the choice of $\w_1$.  Let $k\in K$.  Note that
conjugation by $k^{-1}$ maps $R(\w_1)$ to $R(\w_1k)$.  The
following diagram of isomorphisms commutes.
\begin{equation}\nonumber
\xymatrix{
(R(\w_1)\backslash G_0/K, \omega_{\w_1}) \ar[rr]^{\bu (\w_1(\cdot))} \ar[d]_{conj(k^{-1})} & & (S(\w),\Pi_\w^{-1}) \\
(R(\w_1k)\backslash G_0/K,\omega_{\w_1k}) \ar[urr]_{\bu
(\w_1k(\cdot))} & }
\end{equation}

Our proof of Theorem~\ref{isomtheorem} unfortunately involves a
brutal calculation.  For the convenience of the reader, we
summarize the main steps here.
\begin{list}{$\bullet$}{}
\item First, we introduce an extended operator and a factorization
analogous to the one used in the noncompact case which
we use to compute the inverse of the Poisson tensor on the
symplectic leaf $S(\w)$.
\item Next we compute the derivative
$\Tilde{\bu}_*$ of the map $\Tilde{\bu}\colon R(\w_1)\backslash
G_0/K\to S(\w)$ in our equivariant bundle presentation, and
determine the tangent space to $S(\w)$ in $U\times_K i\p$.  This
is the content of Lemma~\ref{compactcaselemma1}.
\item In Lemma~\ref{lemma42} we produce the expression
\[
\Pi_{\w}^{-1}([u,x]\wedge [u,y])=\langle
\Ad(\w_1g_0)^{-1}\circ\H_{\w}\circ \Ad(\w_1g_0)(x),y\rangle
\]
for the symplectic form $\Pi_\w^{-1}$ on the leaf $S(\w)$.  The
operator $\H_\w$, which arises in the calculation, is precisely
the operator $\H$ in the case $\w=1$.
\item To prove Theorem~\ref{isomtheorem} we must show that
\begin{equation}\label{outline1}
\Pi_{\w}^{-1}(\Tilde{\bu}_*([g_0,x])\wedge\Tilde{\bu}_*([g_0,y]))=\omega_{\w_1}([g_0,x]\wedge
[g_0,y]).
\end{equation}
Writing $\w_1g_0=\bl\ba u$ for the Iwasawa factorization of
$\w_1g_0$ in $G$, the left hand side of (\ref{outline1}) is equal
to the pairing of
\begin{equation}\label{outline2}
\H_\w\circ\Ad(\w_1g_0)(\{\Ad(u^{-1})\circ\H\circ\Ad(u)(x)\}_{\p})
\end{equation}
with
\begin{equation}\label{outline3}
\Ad(\w_1g_0)(\{\Ad(u^{-1})\circ\H\circ\Ad(u)(y)\}_{\p})
\end{equation}
via the Killing form. Parts (a) and (b) of Lemma~\ref{lemma45} and
Lemma \ref{lemma46} are used to simplify the expressions
(\ref{outline2}) and (\ref{outline3}).  Part (c) of Lemma
\ref{lemma45} is used to simplify a later calculation.
\item
Finally, we prove the theorem using the lemmas.
\end{list}


As in the noncompact case (cf. section \ref{The_Noncompact_Case}),
it will be convenient to introduce an extension
$\Tilde{\Omega}(u)$ of $\Omega(u)$ to all of $\u$.  Specifically,
\[
\Tilde{\Omega}(u) = \Ad(u)^{-1}\circ \pr_\u\circ\Ad(u)\circ \Ileft.
\]
A consequence of part (a) of the following lemma is that
\begin{equation}\label{compactextension2}
\Pi_X ([u,\phi]\wedge [u,\psi])=\langle
\Tilde{\Omega}(u)(\phi),\psi \rangle.
\end{equation}

\begin{lemma}\label{lemma40} $\phantom{a}$
\begin{enumerate}
\item[(a)] With respect to the decomposition $\u=\k+i\p$, $\Tilde{\Omega}(u)$ has the form
\begin{equation}\label{compactextension1}
\Tilde{\Omega}(u)=\begin{pmatrix} 1 & b \\ 0 &
\Omega(u)\end{pmatrix}
\end{equation}
for some linear transformation $b\colon i\p\to \k$ depending on
$u$.
\item[(b)] Suppose that $\w_1\in U$ and $g_0\in G_0$ and write $\w_1g_0=\bl\ba u$ for the Iwasawa factorization of $\w_1g_0\in G$.  Then $\Tilde{\Omega}(u)$ can be factored as the composition
\begin{equation}\label{compactextension3}
\xymatrix{ \u \ar[r]^{\Ileft} & \g_0 \ar[r]^{\Ad(g_0)} & \g_0
\ar[r]^{T_{\w_1}(u)} & \u \ar[r]^{\Ad(u)^{-1}} & \u }
\end{equation}
where $T_{\w_1}(u) = \pr_{\u}\circ \Ad((\bl\ba)^{-1}\w_1)$.
\item[(c)] For $X\in\g_0$, $T_{\w_1}(u)(X)\in\u$ is determined by
\begin{equation}\label{compactextension5}
(T_{\w_1}(u)(X))_+ = (((X^{\w_1})_{+})^{(\bl\ba)^{-1}})_{+}
\end{equation}
and
\begin{equation}\label{compactextension6}
(T_{\w_1}(u)(X))_\t=(X^{\w_1}+((X^{\w_1})_{+})^{(\bl\ba)^{-1}})_{\t}.
\end{equation}
\item[(d)] Furthermore, $\ker\Omega(u)=i\{\Ad(\r(\w_1))\}_{\p}$.
\end{enumerate}
\end{lemma}

\begin{proof}
Parts (a), (b), and (d) are immediate consequences of Lemma \ref{lemma22} because  $\Tilde{\Omega}(u)$ is conjugate to $\Tilde{\Omega}_0(u)$ from (\ref{221}).  To be precise,
\[
\Tilde{\Omega}(u)=\Ileft\hspace{-0.5em}\phantom{a}^{-1}\circ\Tilde{\Omega}_0(u)\circ\Ileft .
\]

For part (c), observe that the identity $(X^{(\bl\ba)^{-1}\w_1})_+=(((X^{\w_1})_+)^{(\bl\ba)^{-1}})_+$ is valid because $(\bl\ba)^{-1}\in B^-$.  It then follows that $(T_{\w_1}(X))_+$ is given by (\ref{compactextension5}) and
\begin{eqnarray*}
(T_{\w_1}(u)(X))_\t & = & ((X^{\w_1})^{(\bl\ba)^{-1}})_\t + (((X^{\w_1})_+)^{(\bl\ba)^{-1}})_+ \\
& = & (X^{\w_1})_\t +(((X^{\w_1})_+)^{(\bl\ba)^{-1}})_\t
\end{eqnarray*}
which is equivalent to (\ref{compactextension6}).  This completes the proof.
\end{proof}


\begin{lemma}\label{compactcaselemma1} Assume the hypotheses of Theorem~\ref{isomtheorem}.  Given $u$ representing $uK\in S(\w)$ find $g_0\in G_0$ such that $u=\Tilde{\bu}(g_0)=\bu(\w_1g_0)$.
\begin{enumerate}
\item[(a)] The derivative of the $\tilde {\bu}$-map is given by
\begin{equation}\label{derivativeofutilde}
\begin{matrix}
 G_0\times_K\p & \to & U\times_Ki\p \\
[g_0,y] & \to  & [u,\{\Ad(u^{-1})\circ \H\circ \Ad(u)(-iy)\}_{i\p}].
\end{matrix}
\end{equation}
\item[(b)] The adjoint of the derivative map in (a) is given by
\begin{equation}\nonumber 
\begin{matrix}
U\times_K i\p & \to & G_0\times_K\p \\
 [u,\phi ] &  \to  &  [g_0,\{\Ad(u^{-1})\circ i\H\circ \Ad(u)(\phi )\}_{\p}].
\end{matrix}
\end{equation}\nonumber 
\item[(c)] The tangent space to $S(\w)$ at $uK$ is
\begin{equation}
 \{[u,x]:ix^{g_0}\perp \r(\w_1)\text{ in }\g_0 \}.
\end{equation}
\end{enumerate}
\end{lemma}

\begin{proof} Let $\varepsilon$ denote a small real parameter.  A curve representing the image of $[g_0,y]$ under the derivative of $\Tilde{\bu}$ is
\begin{equation}\label{dercurve}
\varepsilon \mapsto \bu(\w_1g_0e^{\varepsilon y})K.
\end{equation}
Now, $\w_1g_0e^{\varepsilon y}=\bl\ba u e^{\varepsilon y}=\bl\ba
e^{\varepsilon y^u}u$.   Therefore, the linearization of
(\ref{dercurve}) at $\varepsilon = 0$ is given by\[ [u,
\{(\pr_\u(y^u))^{u^{-1}}\}_{i\p}]
\]
since $\bu(\w_1g_0e^{\varepsilon y})$ equals
$uu^{-1}\exp(\varepsilon \pr_{\u}(y^u))u$  to first order.   The
expression in (\ref{derivativeofutilde}) follows using the
commutativity of the left diagram in (\ref{Hisprojuionu}).

Part (b) follows from part (a).  For part (c), we first observe that
\[
T_{uK}(S(\w))=\{[u,\Omega(u)(x)]\colon x\in i\p \}.
\]
The range of $\Omega(u)$ agrees with the orthogonal complement of
its kernel because  $\Omega(u)\in \so(i\p)$.  In
light of part (d) of Lemma \ref{lemma40}, $\ker (\Omega(u))=i\{\Ad(g_0^{-1})(\r(\w_1))\}_\p.$
Let $x\in i\p$ and let $Y\in \r(\w_1)$.  Then
$0=\langle i\{Y^{g_0^{-1}}\}_\p, x\rangle = \langle Y^{g_0^{-1}},
ix \rangle =\langle Y,ix^{g_0}\rangle$
which completes the proof.
\end{proof}

Throughout this section, we will write
\begin{equation}\label{Ad(w)}
\Ad(\w)=\begin{pmatrix} A&0&B\\
0&\hat{w}&0\\
C&0&D\end{pmatrix}
\end{equation}
relative to $\g=\n^{+}+\h+\n^{-}$.  Recall that $\Ad(\w)$ admits
such a presentation because we consider $\w\in N_U(T)$.

\begin{lemma}\label{lemma42} Assume the hypotheses of Theorem \ref{isomtheorem}.  Given $u$ representing $uK\in S(\w)$ find $g_0\in G_0$ such that $u=\Tilde{\bu}(g_0)=\bu(\w_1g_0)$.  Let $[u,x]$ and $[u,y]$ represent tangent vectors to $S(\w)$.

If $\w=\w_1 = 1\in U$ then
\begin{equation}\nonumber 
\Pi_{1}^{-1}([u,x]\wedge [u,y])=\langle \Ad(g_0)^{-1}\circ\H\circ
\Ad(g_0)(x),y\rangle .
\end{equation}
In general
\begin{equation}\nonumber 
\Pi_{\w}^{-1}([u,x]\wedge [u,y])=\langle
\Ad(\w_1g_0)^{-1}\circ\H_{\w}\circ \Ad(\w_1g_0)(x),y\rangle
\end{equation}
where $\H_{\w}:\mathcal D(\H_{\w})\subset \g_0^{\w_1}\to
\g_0^{\w_1}$ is given by
\begin{equation}\nonumber 
\mathcal D(\H_\w)=\{\chi\in\g_0^{\w_1}\colon \chi_-\in
\mathrm{Ran}(1-C\sigma)\}
\end{equation}
\begin{equation}\label{Hwdef}
\H_\w (\chi)= -i\frac {1+C\sigma}{1-C\sigma}\chi_{-}+i\chi_{+}.
\end{equation}
\end{lemma}

\begin{proof} We will prove the general case, as the specific case follows from the fact
that $C=0$ when $\w=1$.  To compute $\Pi_\w^{-1}$ we invert
$\Tilde{\Omega}(u)$ on  the complement of its kernel as
(\ref{compactextension1}) shows that the compression to $i\p$ of
$\Tilde{\Omega}^{-1}(u)$ agrees with $\Omega^{-1}(u)$ on the complement of its kernel.  From (\ref{compactextension2}) and the
factorization of $\Tilde{\Omega}(u)$ in (\ref{compactextension3}),
it follows that
\begin{eqnarray}
\Pi_\w^{-1}([u,x]\wedge [u,y])& = & \langle \Iright\hspace{-0.5em}\phantom{a}^{-1}\circ \Ad(g_0)^{-1}\circ T_{\w_1}^{-1}(u)\circ \Ad(u)(x),y\rangle \nonumber \\
& = & \langle T_{\w_1}^{-1}(u)(x^u),iy^{g_0}\rangle
\label{sympstructlemma1}
\end{eqnarray}
where ``$T_{\w_1}^{-1}(u)(x^u)$'' denotes a solution $X\in \g_0$
to the equation $T_{\w_1}(u)(X)=x^u$.  Such a solution is unique
modulo $\r(\w_1)$.

Write $\chi=x^{\w_1g_0}\in i\g_0^{\w_1}$ and note that
$x^u=\chi^{(\bl\ba)^{-1}}\in \u$.  We seek $X\in\g_0$ such that
$T_{\w_1}(u)(X)=\chi^{(\bl\ba)^{-1}}$.  Using the formulas in
(\ref{compactextension5}) and (\ref{compactextension6}), the
equality of the $\n^+$-components  gives $(X^{\w_1})_+=\chi_+$,
and the equality of the $\h$-components gives that
$\{X^{\w_1}\}_{\t}=\{\chi\}_\t$.  By (c) of Lemma~\ref{compactcaselemma1}, $ix^{g_0}\perp \r(\w_1)$, thus $i\chi
\perp ((\n^-+\a)\cap \g_0^{\w_1})$.  In particular,
$\{i\chi\}_\a=i\{\chi\}_\t=0$, and hence $\{X^{\w_1}\}_\t=0$.

We now know that $X^{\w_1}=L+d+\chi_+$ for some $L\in \n^-$ and
$d\in \g_0^{\w_1}\cap \a$.  The fixed point set of
$\Ad(\w)\circ\sigma$ is $\g_0^{\w_1}$ (this follows from part (a)
of Lemma~\ref{momentummaplemma}).   Thus,
$\Ad(\w)\circ\sigma(X^{\w_1})=X^{\w_1}$ because $X^{\w_1}\in
\g_0^{\w_1}$.  Using the triangular decomposition of $X^{\w_1}$
and the matrix representation of $\Ad(\w)$ in (\ref{Ad(w)}), this
equation implies the following two equations for the $\n^-$ and
$\n^+$-components of $X^{\w_1}$.
\begin{eqnarray}
& & A\sigma(L)=(1-B\sigma)(\chi_+) \label{sympstructlemma2}\\
& & (1-C\sigma)(L)=D\sigma(\chi_+) \label{sympstructlemma3}
\end{eqnarray}
The minus one eigenspace of $\Ad(\w)\circ\sigma$ on $\g$ is
$i\g_0^{\w_1}$ which contains $\chi$.  The equation $\Ad(\w)\circ
\sigma (\chi)=-\chi$ implies the following two equations for the
$\n^-$ and $\n^+$-components of $\chi$.
\begin{eqnarray}
& & B\sigma\chi_+ = -\chi_+-A\sigma\chi_- \label{sympstructlemma4}\\
& & D\sigma\chi_+ = -(1+C\sigma)\chi_- \label{sympstructlemma5}
\end{eqnarray}
Equations (\ref{sympstructlemma3}) and (\ref{sympstructlemma5})
together imply that
\begin{equation}\label{sympstructlemma6}
(1-C\sigma)(L)=-(1+C\sigma)\chi_-
\end{equation}
whereas (\ref{sympstructlemma2}) and (\ref{sympstructlemma4})
together give
\[
A\sigma (L) = 2\chi_+ + A\sigma \chi_-.
\]
The condition that $(1-C\sigma)(L)$ and $A\sigma (L)$ both vanish
is equivalent  to the statement that $L\in \n^-$ belongs to
$\g_0^{\w_1}$.  Thus (\ref{sympstructlemma6}) implies that
\begin{equation}\nonumber
L = -\frac{1+C\sigma}{1-C\sigma}\chi_- \text{ modulo } \n^- \cap
\g_0^{\w_1}.
\end{equation}

Note that $d\in \a\cap\g_0^{\w_1}$ so
\[
X^{\w_1} = -\frac{1+C\sigma}{1-C\sigma}\chi_-+\chi_+\text{ modulo }(\n+\a)\cap\g_0^{\w_1}=\r(\w_1)^{\w_1}.
\]
Therefore $X=\Ad(\w_1^{-1})\circ\H_\w(-i\chi)$ modulo $\r(\w_1)$ where $\H_\w$ is as in (\ref{Hwdef}).  Substituting this for
$T_{\w_1}^{-1}(u)(x^u)$ in  (\ref{sympstructlemma1}) completes the proof
of Lemma~\ref{lemma42}.
\end{proof}

In the hypothesis of Theorem~\ref{isomtheorem}, we consider
$\w=\w_1\w_1^{-\Theta}$.  By part (a) of Lemma~\ref{momentummaplemma}, $\Ad(\w)\circ\Theta$ is a complex linear
involution of $\g$ which commutes with the Cartan involution
fixing $\u$.  The composition of $\Ad(\w)\circ \Theta$ with the
Cartan involution is the complex anti-linear involution
$\Ad(\w)\circ\sigma$ and its fixed point set is the real form
$\g_0^{\w_1}$. Given $Z\in \g$, we denote its orthogonal
projection to $\g_0^{\w_1}$ by
\begin{equation}\nonumber 
\{Z\}_{\g_0^{\w_1}}=\textstyle{\frac{1}{2}}(Z+\Ad(\w)\circ\sigma(Z)).
\end{equation}
Note that $\Ad(\w_1)$ intertwines the orthogonal projections to
$\g_0$ and $\g_0^{\w_1}$, i.e.,
\begin{equation*}
\Ad(\w_1)(\{Z\}_{\g_0})=\{\Ad(\w_1)(Z)\}_{\g_0^{\w_1}}\quad
\forall \,Z\in \g.
\end{equation*}

\begin{lemma}\label{lemma45}
 Assume the hypotheses of Theorem~\ref{isomtheorem}.
 Let $p_-$, $p_0$, and $p_+$ denote the projections corresponding to
the triangular decomposition $\g=\n^-+\h+\n^+$, write
$\w_1g_0=\bl\ba u$ for the Iwasawa factorization of $\w_1g_0$ in
$G$, and let $x\in\p$. Then
\begin{enumerate}
\item[(a)] $\Ad(\bl\ba)\circ\H\circ\Ad(\bl\ba)^{-1}-\H=i\mathcal Z$
where $\mathcal Z\colon \g\to \n^-+\h$ is the operator
\begin{eqnarray}
\mathcal Z & = & -p_0\circ \Ad(\bl\ba)^{-1}\circ p_{+}+2p_{-}\circ \Ad(\bl\ba)\circ p_{+}\circ \Ad(\bl\ba)^{-1}\circ p_{+}+ \label{defofcalZ}\\
& & p_{-}\circ\Ad(\bl\ba)\circ p_0\circ \Ad(\bl\ba)^{-1}\circ
p_{+}+p_{-}\circ \Ad(\bl\ba)\circ p_0, \nonumber
\end{eqnarray}
\item[(b)] the value of $\mathcal Z$ on $x^{\w_1g_0}$ is
\begin{equation}\label{valueofZ}
\mathcal
Z(x^{\w_1g_0})=-((x^{u})_0-(x^{\w_1g_0})_0)+p_-\circ\Ad(\bl\ba)((x^{u})_0+2(x^{u})_{+}),
\end{equation}
\item[(c)] and for $\chi\in \b^-$
\begin{equation}\nonumber 
\H_\w(\{\chi\}_{i\g_0^{\w_1}}) = \{-i\,p_-(\chi)\}_{i\g_0^{\w_1}}.
\end{equation}

\end{enumerate}
\end{lemma}

\begin{proof}
Relative to the triangular decomposition, written in the order
$\n^{+}+\h+\n^{-}$, $\Ad(\bl\ba)$, $\H$, and $\Ad(\bl\ba)^{-1}$
are represented as $3\times 3$ matrices,
\[
\left(\begin{matrix} \mu&0&0\\
\mu'&1&0\\
\lambda&\nu'&\nu\end{matrix} \right),\quad\left(\begin{matrix} i&0&0\\
0&0&0\\
0&0&-i\end{matrix} \right),\text{ and }\quad\left(\begin{matrix} M&0&0\\
M'&1&0\\
L&N'&N\end{matrix} \right),
\]
respectively, where $\mu =p_{+}\circ \Ad(\bl\ba)\circ p_{+}$, etc.
Note that
$$\mu'M+M',\quad\lambda M+\nu'M'+\nu L,\quad and\quad\nu'
+\nu N'$$ all vanish.  Then
$$\Ad(\bl\ba)\circ \H\circ \Ad(\bl\ba)^{-1}-\H
=\left(\begin{matrix} 0&0&0\\
i\mu'M&0&0\\
i(\lambda M-\nu L)&-i\nu N'&0\end{matrix} \right)$$

$$=\left(\begin{matrix} 0&0&0\\
i\mu'M&0&0\\
i(2\lambda M+\nu'M')&-i\nu N'&0\end{matrix}
\right)=\left(\begin{matrix}
0&0&0\\
-iM'&0&0\\
i(2\lambda M+\nu'M')&i\nu'&0\end{matrix} \right).$$
This last expression is equivalent to (\ref{defofcalZ}).  The proof of part
(a) is complete.

For part (b), we evaluate $\mathcal Z$ from part (a) on
$x^{\w_1g_0}$ and obtain a sum of four terms.
\begin{eqnarray}
&  &  -p_0\circ \Ad(\bl\ba)^{-1}((x^{\w_1g_0})_{+})+2p_{-} \circ
\Ad(\bl\ba)\circ p_{+}\circ \Ad(\bl\ba)^{-1}((x^{\w_1g_0})_{+}) \label{RHSofZvalue}\\
& & + p_{-}\circ\Ad(\bl\ba)\circ p_0\circ
\Ad(\bl\ba)^{-1}((x^{\w_1g_0})_{+})+p_{-}\circ \Ad(\bl\ba)((x^{\w_
1g_0})_0)\nonumber
\end{eqnarray}
The first three terms involve the expression
$\Ad(\bl\ba)^{-1}((x^{\w_1g_0})_+)$ which equals
\[
x^u-\Ad(\bl\ba)^{-1}((x^{\w_1g_0})_0)-\Ad(\bl\ba)^{-1}((x^{\w_1g_0})_-).
\]
The zero mode of this expression is $(x^u)_0-(x^{\w_1g_0})_0$ and
the $\n^+$ projection is $(x^u)_+$ since $\bl\ba\in B^-$.
Inserting these computations into the right hand side of
(\ref{RHSofZvalue}) gives that
\begin{eqnarray*}
\mathcal Z(x^{\w_1g_0}) & = & -((x^u)_0-(x^{\w_1g_0})_0)+2p_-\circ\Ad(\bl\ba)((x^u)_{+})\\
&&+p_-\circ\Ad(\bl\ba)((x^u)_0-(x^{\w_1g_0})_0)+
p_-\circ\Ad(\bl\ba)((x^{\w_1g_0})_0)\\
&=&
-(x^u)_0+(x^{\w_1g_0})_0+p_-\circ\Ad(\bl\ba)((x^u)_0+2(x^u)_{+}).
\end{eqnarray*}

For part (c), let $\chi\in \b^-$ and write $\chi=\chi_-+\chi_0$ for
its triangular decomposition.  Using the matrix (\ref{Ad(w)})
representing $\Ad(\w)$ and the fact that $\sigma(\n^-)=\n^+$, one
computes that
\begin{eqnarray*}
\{\chi\}_{i\g_0^{\w_1}}& =& \frac{1}{2}(\chi -\Ad(\w)\sigma (\chi)) \\
& = & \frac{1}{2}( (1-C\sigma)(\chi_-)+\hat{w}\sigma (\chi_0) +
A\sigma (\chi_-)).
\end{eqnarray*}
Observe that $\{\chi\}_{i\g_0^{\w_1}}$ is in the domain of the
operator $\H_\w$ from Lemma \ref{lemma42}.  From the definition of
that operator (\ref{Hwdef}) we see that $\H_\w$ multiplies the
upper triangular part by $i$ and the lower triangular part by the
operator $-i(1+C\sigma)(1-C\sigma)^{-1}$ and kills the zero mode.
Hence,
\begin{eqnarray*}
\H_\w(\{\chi\}_{i\g_0^{\w_1}}) & = & \frac{1}{2}(-i(1+C\sigma)(\chi_-)+iA\sigma(\chi_-)) \\
& = & \frac{1}{2}((1-C\sigma)(-i\chi_-)+A\sigma(-i\chi_-)) \\
& = & \{-i\chi_-\}_{i\g_0^{\w_1}}.
\end{eqnarray*}
This completes the proof of Lemma~\ref{lemma45}.
\end{proof}

\begin{lemma}\label{lemma46}
Assume the hypotheses of Theorem~\ref{isomtheorem} and suppose
that $x\in \p$.  Then
\begin{equation*}
\Ad(\w_1g_0)(\{ \Ad(u^{-1})\circ \H\circ
\Ad(u)(x)\}_{\p})=-i\{\Ad(\bl\ba)\circ\pr_{\n^-+\a}(x^u)\}_{i\g_0^{\w_1}}.
\end{equation*}
\end{lemma}
\begin{proof}
Since $x\in\p\subset i\u$, and $\H$ and $\Ad(u)$ preserve $i\u$,
it follows that the orthogonal projection to $\p$ of
$\Ad(u)^{-1}\circ\H\circ\Ad (u)(x)$ agrees with its orthogonal
projection to $\g_0$.  We compute:
\begin{eqnarray}
& & \Ad(\w_1g_0)(\{\Ad(u^{-1})\circ \H\circ \Ad(u)(x)\}_{\p})\nonumber \\
&= & \Ad(\w_1g_0)(\{\Ad(u^{-1})\circ \H\circ \Ad(u)(x)\}_{\g_0})\nonumber \\
&=& \{\Ad(\w_1g_0)\circ\Ad(u^{-1})\circ \H\circ
\Ad(u)\circ\Ad(\w_1g_0)^{-1}(x^{\w_1g_0})\}_{\g_0^{\w_1}}\label{lemma46_eqA}
\end{eqnarray}
where in the last line we used that $\Ad(g_0)$ commutes with the
orthogonal projection to $\g_0$ and $\Ad(\w_1)$ intertwines the
projections to $\g_0$ and $\g_0^{\w_1}$.  Using the factorization
$\w_1g_0=\bl\ba u$ we obtain that (\ref{lemma46_eqA}) is equal to
\begin{eqnarray}
&=&\{(\Ad(\bl\ba)\circ\H\circ\Ad(\bl\ba)^{-1}(x^{\w_1g_0})\}_{\g_0^{\w_1}}\nonumber \\
&=&\{(\Ad(\bl\ba)\circ\H\circ\Ad(\bl\ba)^{-1}-\H)(x^{\w_1g_0})\}_{\g_0^{\w_1}}+\{\H(x^{\w_1g_0})\}_{\g_0^{\w_1}}\nonumber \\
&=&\{i\mathcal
Z(x^{\w_1g_0})\}_{\g_0^{\w_1}}+\{\H(x^{\w_1g_0})\}_{\g_0^{\w_1}}
\label{lemma46_eqB}
\end{eqnarray}
where $\mathcal Z$ is the operator from part (a) of Lemma~\ref{lemma45}. Examine the value of $i\mathcal Z(x^{\w_1g_0})$
using (\ref{valueofZ}).  The zero mode is a sum of two terms,
namely $-(ix^u)_0$ and $(ix^{\w_1g_0})_0$.  The latter of these
is in $i\g_0^{\w_1}$ which is the kernel of the orthogonal
projection to $\g_0^{\w_1}$.  The former is in $\t$, and thus
$\{-(ix^u)_0\}_{\g_0^{\w_1}}=-\{ix^u\}_{\t_w}$.  Thus
(\ref{lemma46_eqB}) becomes
\begin{equation}\label{lemma46_eqC}
\{-ix^u\}_{\t_w}+\{p_-\circ\Ad(\bl\ba)(i(x^{u})_0+2i(x^{u})_{+})\}_{\g_0^{\w_1}}+\{\H(x^{\w_1g_0})\}_{\g_0^{\w_1}}.
\end{equation}
Now we assert that for $\chi\in \g_0^{\w_1}$,
$\{\H(\chi)\}_{\g_0^{\w_1}}=\{-2ip_-(\chi)\}_{\g_0^{\w_1}}$.  To
see this, note that $\Ad(\w)\circ\sigma (\chi)=\chi$ implies that
$\chi_0\in \t_w$, and $\chi_-$ and $\chi_+$ satisfy the equations
\begin{eqnarray*}
A\sigma(\chi_-) & = & (1-B\sigma)(\chi_+) \\
\text{and }(1-C\sigma)(\chi_-) & = & D\sigma(\chi_+).
\end{eqnarray*}
Computing the orthogonal projection of
$\H(\chi)=-i\chi_-+i\chi_+$, keeping in mind the above relations,
and using the complex anti-linearity of $\sigma$ establishes the
assertion.

We now apply this assertion to
\[
p_-(x^{\w_1g_0})=p_-(x^{\bl\ba
u})=p_-\circ\Ad(\bl\ba)((x^u)_-+(x^u)_0+(x^u)_+)
\]
obtaining that $\{\H(x^{\w_1g_0})\}_{\g_0^{\w_1}}$ equals
\begin{equation}\label{lemma46_eqD}
 \{-2ip_-(x^{\w_1g_0})\}_{\g_0^{\w_1}}=\{p_-\circ\Ad(\bl\ba)(-2i(x^u)_--2i(x^u)_0-2i(x^u)_+)\}_{\g_0^{\w_1}}.
\end{equation}
Replacing the third term of (\ref{lemma46_eqC}) with the right
hand side of (\ref{lemma46_eqD}) and combining the terms gives
that (\ref{lemma46_eqC}) is equal to
\begin{equation}\label{lemma46_eqE}
 -\{ix^u\}_{\t_w}+\{p_-\circ\Ad(\bl\ba)(-2i(x^u)_--i(x^u)_0)\}_{\g_0^{\w_1}}.
\end{equation}
Observe that $x^u\in i\u$ and thus
$\pr_{\n^-+\a}(x^u)=2(x^u)_-+(x^u)_0$.  Hence (\ref{lemma46_eqE})
may be rewritten as
\begin{eqnarray}
& = & -\{ix^u\}_{\t_w}+\{-i\,p_-\circ\Ad(\bl\ba)\circ\pr_{\n^-+\a}(x^u)\}_{\g_0^{\w_1}} \nonumber \\
& = &
-i\{x^u\}_{i\t_w}-i\{\,p_-\circ\Ad(\bl\ba)\circ\pr_{\n^-+\a}(x^u)\}_{i\g_0^{\w_1}}.
\label{lemma46_eqF}
\end{eqnarray}
Notice that
\begin{eqnarray*}
i\{p_-\circ\Ad(\bl\ba)\circ\pr_{\n^-+\a}(x^u)\}_{i\g_0^{\w_1}} & = & i\{\Ad(\bl\ba)\circ\pr_{\n^-+\a}(x^u) -(x^u)_0 \}_{i\g_0^{\w_1}} \\
& = & i\{\Ad(\bl\ba)\circ\pr_{\n^-+\a}(x^u)\}_{i\g_0^{\w_1}}
-i\{x^u\}_{i\t_w}
\end{eqnarray*}
because $\bl\ba\in B^-$.  Therefore, (\ref{lemma46_eqF}) equals
\[
-i\{\Ad(\bl\ba)\circ\pr_{\n^-+\a}(x^u)\}_{i\g_0^{\w_1}}.
\]
This proves Lemma~\ref{lemma46}.
\end{proof}

\begin{proof}[Proof of Theorem~\ref{isomtheorem}]
We consider the symplectic leaf $S(\w)$ and fix a choice of
$\w_1\in U$ such that $\w_1\w_1^{-\Theta}=\w\in N_U(T)$.  Write
$w$ for the element of the Weyl group represented by $\w$.  Our
goal is to show that under the map $\Tilde{\bu}\colon G_0/K\to
U/K: g_0K\mapsto uK$, where $\w_1g_0=\bl\ba u$, the symplectic
form $\Pi^{-1}_\w$ pulls back to the global two-form
$\omega_{\w_1}$ from (\ref{defofomega}), i.e., we need to show
that
\begin{equation}\nonumber
\Pi^{-1}_\w(\Tilde{\bu}_*[g_0,x]\wedge\Tilde{\bu}_*[g_0,y])=\langle
\Ad(u^{-1})\circ \H\circ \Ad(u)(x),y\rangle.
\end{equation}
By part (a) of Lemma~\ref{compactcaselemma1},
\begin{eqnarray*}
\Tilde{\bu}_*[g_0,x] & = &[u,\{\Ad(u^{-1})\circ\H\circ\Ad(u)(-ix)\}_{i\p}]\\
& = & [u,-i\{\Ad(u^{-1})\circ\H\circ\Ad(u)(x)\}_{\p}].
\end{eqnarray*}
We will write
\begin{eqnarray*}
X & = & \Ad(\w_1g_0)(-i\{\Ad(u^{-1})\circ\H\circ\Ad(u)(x)\}_{\p})\\
\text{and }Y & = &
\Ad(\w_1g_0)(-i\{\Ad(u^{-1})\circ\H\circ\Ad(u)(y)\}_{\p}).
\end{eqnarray*}
Then, by Lemma~\ref{lemma42},
$\Pi^{-1}_\w(\Tilde{\bu}_*[g_0,x]\wedge\Tilde{\bu}_*[g_0,y])=\langle
\H_\w(X),Y\rangle $. Using Lemma~\ref{lemma46}, the expressions
for $X$ and $Y$ may be simplified to
\[
X=-\{\Ad(\bl\ba)\circ\pr_{\n^-+\a}(x^u)\}_{i\g_0^{\w_1}}
\]
and
\[
Y=-\{\Ad(\bl\ba)\circ\pr_{\n^-+\a}(y^u)\}_{i\g_0^{\w_1}}.
\]
Notice that $X$ is the projection to $i\g_0^{\w_1}$ of an element
of $\b^-$.  Thus, by part (c) of Lemma~\ref{lemma45},
\[
\H_\w(X)=\{i
\,p_-\circ\Ad(\bl\ba)\circ\pr_{\n^-+\a}(x^u)\}_{i\g_0^{\w_1}}.
\]
The subalgebra $\g_0^{\w_1}$ is a real form of $\g$, and the
Killing form is real valued on the real subspace $i\g_0^{\w_1}$.
Therefore,
\begin{eqnarray}
\langle \H_\w(X),Y\rangle & = & \langle \{i \,p_-\circ \Ad(\bl\ba)\circ\pr_{\n^-+\a}(x^u)\}_{i\g_0^{\w_1}}, -\{\Ad(\bl\ba)\circ\pr_{\n^-+\a}(y^u)\}_{i\g_0^{\w_1}}\rangle \nonumber \\
& = & \Re\langle i \,p_-\circ\Ad(\bl\ba)\circ\pr_{\n^-+\a}(x^u), -\{\Ad(\bl\ba)\circ\pr_{\n^-+\a}(y^u)\}_{i\g_0^{\w_1}}\rangle \nonumber \\
& = & \Im\langle p_-\circ\Ad(\bl\ba)\circ\pr_{\n^-+\a}(x^u),
\{\Ad(\bl\ba)\circ\pr_{\n^-+\a}(y^u)\}_{i\g_0^{\w_1}}\rangle
.\label{myproof1}
\end{eqnarray}
Applying the definition of the orthogonal projection to
$i\g_0^{\w_1}$, and expanding (\ref{myproof1}) using bilinearity
of the Killing form, we obtain a sum of two terms.  The first one
vanishes because $\Im\langle \n^-,\n^-+\a\rangle =0$.  So
(\ref{myproof1}) is equivalent to
\begin{eqnarray}
& = & \frac{1}{2}\Im\langle p_-\circ \Ad(\bl\ba)\circ \pr_{\n^-+\a}(x^u),\Ad(\bl\ba)\circ\pr_{\n^-+\a}(y^u) \rangle \nonumber \\
&  & -\frac{1}{2}\Im\langle p_-\circ\Ad(\bl\ba)\circ\pr_{\n^-+\a}(x^u),\Ad(\w)\circ\sigma\circ\Ad(\bl\ba)\circ\pr_{\n^-+\a}(y^u) \rangle \nonumber \\
& = & -\frac{1}{2}\Im\langle
p_-\circ\Ad(\bl\ba)\circ\pr_{\n^-+\a}(x^u),\Ad(\w)\circ\sigma\circ\Ad(\bl\ba)\circ\pr_{\n^-+\a}(y^u)
\rangle . \label{myproof2}
\end{eqnarray}

We can replace $p_-\circ\Ad(\bl\ba)\circ\pr_{\n^-+\a}(x^u)$ by
$\Ad(\bl\ba)\circ\pr_{\n^-+\a}(x^u)-(x^u)_0$ in the left hand
factor of (\ref{myproof2}) since $\bl\ba\in B^-$ and expand again
continuing our calculation:
\begin{eqnarray}
& = & -\frac{1}{2}\Im\langle \Ad(\bl\ba)\circ\pr_{\n^-+\a}(x^u),\Ad(\w)\circ\sigma\circ\Ad(\bl\ba)\circ\pr_{\n^-+\a}(y^u) \rangle \label{myproof3}\\
&  & +\frac{1}{2}\Im\langle
(x^u)_0,\Ad(\w)\circ\sigma\circ\Ad(\bl\ba)\circ\pr_{\n^-+\a}(y^u)
\rangle . \nonumber
\end{eqnarray}
The second term vanishes because the zero mode of the right hand
factor is in $\a$ and thus its pairing with $(x^u)_0\in \a$ is
real, so the imaginary part is zero.  From the factorization
$\w_1g_0=\bl\ba u$, and the equation $\w_1\w_1^{-\Theta}=\w$, we
have that $(\bl\ba)^{-1}\w(\bl\ba)^\sigma=uu^{-\Theta}$.  Thus,
(\ref{myproof3}) is equivalent to
\begin{eqnarray}
& = & -\frac{1}{2}\Im\langle \Ad(\bl\ba)\circ\pr_{\n^-+\a}(x^u),\Ad(\w)\circ\sigma\circ\Ad(\bl\ba)\circ\pr_{\n^-+\a}(y^u) \rangle \nonumber \\
& = & -\frac{1}{2}\Im\langle
\pr_{\n^-+\a}(x^u),\Ad(uu^{-\Theta})\circ\sigma\circ\pr_{\n^-+\a}(y^u)
\rangle. \label{myproof4}
\end{eqnarray}
Now recall from Proposition~\ref{prop12} that
$\pr_\u(\cdot)=\H(-i\cdot)=-i\H(\cdot)$ on $i\u$.  Since $x^u\in i\u$,
$\pr_{\n^-+\a}(x^u)=x^u-\pr_\u(x^u)=x^u+i\H(x^u)$, and similarly
$\pr_{\n^-+\a}(y^u)=y^u+i\H(y^u)$. Making these replacements in
(\ref{myproof4}), we obtain
\begin{eqnarray}
& =  & -\frac{1}{2}\Im\langle x^u+i\H(x^u),\Ad(uu^{-\Theta})\circ\sigma(y^u+i\H(y^u)) \rangle \label{myproof5} \\
& = & -\frac{1}{2}\Im\langle x^u ,\Ad(uu^{-\Theta})\circ\sigma(i\H(y^u)) \rangle \label{myproof6} \\
& & \quad -\frac{1}{2}\Im\langle i\H(x^u),
\Ad(uu^{-\Theta})\circ\sigma(y^u)\rangle \nonumber
\end{eqnarray}
using that $\Im\langle i\u,i\u\rangle = 0$.  Since $\sigma$ agrees
with $\Theta$ on $U$ and fixes $\p\subset\g_0$, it follows that
$\Ad(uu^{-\Theta})\circ\sigma$ fixes $y^u$.
Therefore (\ref{myproof6})
\begin{eqnarray}
& = & -\frac{1}{2}\Im\langle \Ad(u^\Theta)(x) ,\sigma(i\H(y^u)) \rangle -\frac{1}{2}\Im\langle i\H(x^u), y^u\rangle \nonumber \\
& = & -\frac{1}{2}\Im\langle \sigma (x^u) ,\sigma(i\H(y^u))
\rangle -\frac{1}{2}\Im\langle i\H(x^u), y^u\rangle
\label{myproof7}
\end{eqnarray}
where (\ref{myproof7}) is obtained from the previous line using
the facts that $\sigma\circ\Ad(u)=\Ad(u^\Theta)\circ \sigma$ and
$\p$ is fixed by $\sigma$.  Finally, we have
that (\ref{myproof7})
\begin{eqnarray*}
& = & -\frac{1}{2}\langle x^u ,\H(y^u) \rangle +\frac{1}{2}\langle \H(x^u), y^u\rangle \\
& = &\langle \Ad(u^{-1})\circ\H\circ\Ad(u)(x),y\rangle
\end{eqnarray*}
because of the skew-symmetry of $\H$. The proof of Theorem~\ref{isomtheorem} is now complete.
\end{proof}

\section{The Group Case}\label{The_Group_Case}

Let $K$ be a simply connected compact Lie group.  With respect to
the invariant metric induced by the Killing form, $K$ may be
viewed as a compact symmetric space $X$.  In this case, the
diagram in (\ref{diagram_of_groups}) specializes to
\begin{equation}\nonumber
\xymatrix{
 & G=K^\C \times K^\C &  \\
K^\C \simeq G_0 \ar[ur] & & U= K\times K \ar[ul]\\
 & \Delta(K)\ar[ul] \ar[ur] &
}
\end{equation}
where $\Delta(K)=\{(k,k)\colon k\in K\}$ and
$G_0=\{g_0=(g,g^{-*})\colon g\in K^\C\}$.  The involution $\Theta$
in this case is the outer automorphism $\Theta
((g_1,g_2))=(g_2,g_1)$. Also
\[
X_0=G_0/\Delta(K)\simeq K^{\C}/K,\text{ and } X=U/\Delta(K)\simeq
K,
\]
where the latter isometry is $(k_1,k_2)\Delta(K)\mapsto
k=k_1k_2^{-1}$.

To distinguish between $\g$ and $\k^{\C}$, we will adopt the
(admittedly cumbersome) convention of denoting structures
associated with $\k^{\C}$ using superchecks.

We fix a triangular decomposition
\begin{equation}\label{groupcase1}
\check {\g}=\k^{\C}=\check{\n}^{-} + \check{\h} + \check {\n}^{+}.
\end{equation}
This induces a $\Theta$-stable triangular decomposition for $\g$
\begin{equation}\label{groupcase2}
\g =
\underbrace{(\check{\n}^-\times\check{\n}^-)}_{\n^-}+\underbrace{(\check{\h}\times\check{\h})}_{\h}+\underbrace{(\check{\n}^+\times\check{\n}^+)}_{\n^+}
\end{equation}
Let $\check {\a}=\check {\h}_{\R}$ and $\check {\t}=i\check {\a}$.
Then
\[
\t_0=\{(x,x):x\in\check {\t}\},\quad\text{\rm and}\quad
\a_0=\{(y,-y):y\in\check {\a}\}.
\]
The standard Poisson Lie group structure on $U=K\times K$ induced
by the decomposition in (\ref{groupcase2}) is then the product
Poisson Lie group structure for the standard Poisson Lie group
structure on $K$ induced by the decomposition (\ref{groupcase1}).

Let us denote the Poisson Lie group structure on $K$ by $\pi_K$
and the Evens-Lu homogeneous Poisson structure on $X=K$ by
$\Pi_K$.  The identification of $\k$ with its dual via the Killing
form allows us to view the Hilbert transform $\check{\H}$
associated to (\ref{groupcase1}) as an element of $\k\wedge\k$. As
a bivector field
\[
\pi_K = \check{\H}^r-\check{\H}^l
\]
where $\check{\H}^r$ (resp. $\check{\H}^l$) denotes the right
(resp. left) invariant bivector field on $K$ generated by
$\check{H}$, whereas $\Pi_K=\check{\H}^r+\check{\H}^l$ (see
section 5 of \cite{Caine}).

\begin{theorem}\label{bruhatbirkhoff}
Let $\w_0\in N_K(T)$ be a representative for the longest element
of the Weyl group.  The map $L_{\w_0}\colon K\to K$
\begin{equation}\label{groupcase3}
K\ni k\mapsto \w_0k\in K
\end{equation}
is a Poisson diffeomorphism carrying the Poisson Lie group
structure $\pi_K$ onto (the negative of) the Evens-Lu $(K\times
K,\pi_K\oplus\pi_K)$-homogeneous Poisson structure $\Pi_X$ on
$X=K$.
\end{theorem}

\begin{proof}
Identify the dual of $\k$ with $\k$ using the Killing form, and
use right translation to trivialize $T^*K$ as $K\times \k$.  From
section 5 of \cite{Caine}, we have that
\begin{equation}\nonumber 
\pi_K((k,\phi),(k,\psi))=\langle
(\check{\H}-\Ad(k)\circ\check{\H}\circ\Ad(k)^{-1})(\phi),\psi\rangle
\end{equation}
and
\begin{equation}\nonumber 
\Pi_K((k,\phi),(k,\psi))=\langle
(\check{\H}+\Ad(k)\circ\check{\H}\circ\Ad(k)^{-1})(\phi),\psi\rangle
\end{equation}
for each $(k,\phi)$ and $(k,\psi)$ representing cotangent vectors
at $k\in K$.  With the tangent bundle to $K$ trivialized as
$K\times \k$ using right translation, the derivative of
(\ref{groupcase3}) is $(k,X)\mapsto (\w_0k,\Ad(\w_0)(X))$ and the
transpose map is
\begin{equation}\nonumber
(\w_0k,\phi)\mapsto (k,\Ad(\w_0)^{-1}(\phi)).
\end{equation}
Then $L_{\w_0*}\pi_K ((\w_0k,\phi),(\w_0k,\psi)) $ is
\begin{eqnarray}
 & = & \langle \check{\H}\circ \Ad(\w_0)^{-1}(\phi), \Ad(\w_0)^{-1}(\psi)\rangle \nonumber \\
&  & - \langle \Ad(k)\circ\check{\H}\circ \Ad(k)^{-1}\circ\Ad(\w_0)^{-1}(\phi), \Ad(\w_0)^{-1}(\psi)\rangle \nonumber \\
& = & \langle \Ad(\w_0)\circ\check{\H}\circ\Ad(\w_0)^{-1}(\phi),\psi\rangle \nonumber \\
& & -\langle
\Ad(\w_0k)\circ\check{\H}\circ\Ad(\w_0k)^{-1}(\phi),\psi\rangle.
\label{groupcase6}
\end{eqnarray}
The operator $\check{\H}$ is conjugated to $-\check{\H}$ by
$\Ad(\w_0)$ as conjugation by $\w_0$ interchanges $\check{\n}^-$
and $\check{\n}^+$.  Thus (\ref{groupcase6}) becomes
\begin{eqnarray*}
& & -\langle (\check{\H}+\Ad(\w_0k)\circ\check{\H}\circ\Ad(\w_0k)^{-1})(\phi),\psi\rangle \\
& = & -\Pi_K ((\w_0k,\phi),(\w_0k,\psi)).
\end{eqnarray*}
This completes the proof.
\end{proof}

The symplectic leaves of $\Pi_K$ foliate the strata of the Birkoff
decomposition of $K$ induced by (\ref{groupcase1}).  The top
stratum, $\Sigma_1^K$, consists of those elements admitting a
unique triangular factorization $k=lmau$ where $l\in \check{N}^-$,
$m\in \exp(\check{t})=\check{T}$, $a\in \exp(\check{a})=\check{A}$
and $u\in\check{N}^+$.  This stratum is an open dense subset of
$K$.

The symplectic leaf through the identity, $S(1)$, consists of
those elements whose factorization has $m=1$. In the remainder of
this section, we will focus on this one leaf. We will generally
identify this leaf with $\check{N}^-$, using $l$ as a global
coordinate. The Hamiltonian action of $T_0$ on this leaf is
isomorphic to the conjugation action of $\check {T}$ on
$l\in\check{N}^-$.

The isomorphism in Theorem~\ref{isomtheorem} is given by
$$\check {A}\backslash K^{\C }/K\to A\backslash G_0/\Delta (K)
\to S(1)$$
$$\check {A}gK\mapsto A(g,g^{-*})\Delta (K)\mapsto k=a_1^{-1}l_1^{-1}l_2^{
-*}a_2^{-1},$$ in terms of the Iwasawa decompositions
$$g=l_1a_1k_1,\quad g^{-*}=l_2a_2k_2.$$
We can clearly take $g=l_1=l^{-1}$, implying that
$$a=a_2^{-1},\quad u=a_2l_2^{-*}a_2^{-1}.$$
From the noncompact perspective, $l$ is essentially a horocycle
coordinate, and from the compact perspective, $l$ is a standard
affine coordinate for the flag space $K/T$.

In the case $K=SU(2)$, this coordinate is given explicitly by
\begin{equation}\nonumber 
l=\left(\begin{matrix} 1&0\\
\zeta&1\end{matrix} \right)\leftrightarrow k(\zeta
)=\left(\begin{matrix}
1&0\\
\zeta&1\end{matrix} \right)\left(\begin{matrix} a&0\\
0&a^{-1}\end{matrix} \right)\left(\begin{matrix} 1&-\bar{\zeta}\\
0&1\end{matrix} \right),\end{equation} where
$a=(1+\vert\zeta\vert^2)^{-1/2}$.

The following theorem is a reformulation of results in \cite{Lu}
on the standard Poisson structure. This reformulation is of
importance in connection with infinite dimensional generalizations
(see \cite{pickrell2}). We denote the symplectic form on $\check
{N}^{-}$ simply by $\omega$ (in our earlier notation this is
$\omega_1$, from the noncompact point of view, and $\Pi_1^{-1}$, from
the compact point of view). We assume that we are given a Serre
presentation compatible with the triangular decomposition of
$\k^{\C}$. Given a simple positive root $\gamma$, we let
$i_{\gamma}:SU(2) \mapsto K$ denote the corresponding root
subgroup inclusion, and
$$r_{\gamma}=i_{\gamma}(\left(\begin{matrix} 0&i\\i&0\end{matrix}\right)),$$
a fixed representative for the corresponding Weyl group
reflection.

\begin{theorem}\label{Lutheorem} Fix $w\in W$.
\begin{enumerate}
\item[(a)]  The submanifold ($\check {N}^{-}\cap w^{-1}\check
{N}^{+}w)\subset \check {N}^{-}$ is $\check {T}$-invariant and
symplectic.

Fix a representative $\w$ for $w$ with minimal factorization
$\w=r_n..r_1$, in terms of simple reflections $r_j=r_{\gamma_j}$
corresponding to simple positive roots $\gamma_j$.  Let
$w_j=r_j..r_1$.

\item[(b)]  The map
$$\C^n\to N^{-}\cap w^{-1}N^{
+}w:\zeta =(\zeta_n,\dots,\zeta_1)\to l(\zeta)$$ where
$$w_{n-1}^{-1}i_{\gamma_
n}(k(\zeta_n))w_{n-1}\dots w_1^{-1}i_{\gamma_2}(k(\zeta_2))w_1i_{\gamma_1}(k(\zeta_1))=l(\zeta )au,$$
is a diffeomorphism.

\item[(c)] In these coordinates the restriction of $\omega$ is given by
$$\omega\vert_{N^{-}\cap w^{-1}N^{+}w}=\sum_{j=1}^n\frac i{\langle\gamma_
j,\gamma_j\rangle}\frac 1{(1+\vert\zeta_j\vert^2)}d\zeta_j\wedge
d\bar{\zeta}_j,$$ the momentum map is the restriction of
$-\langle\frac i2log(a ),\cdot\rangle$, where
$$a(k(\zeta ))=\prod_{j=1}^n (1+\vert\zeta_j\vert^2)^{-\frac 12 w_{j-1}^{-1}
h_{\gamma_j}w_{j-1}},$$
 and Haar measure (unique up to a constant) is given by

$$d\lambda_{N^{-}\cap w^{-1}N^{+}w}(l)=\prod_{j=1}^n (1+\vert\zeta_j\vert^2)^{\check{\delta}(w_{j-1}^{-
1}h_{\gamma_j}w_{j-1})-1},$$  where $\check{\delta
}=\sum\Lambda_j$, the sum of the dominant integral functionals for
$\check {g}$, relative to (\ref{groupcase1}).

\item[(d)] Let $C_{\w}$ denote the symplectic leaf through $\w$, with
respect to $\Pi_K$, with the negative of the induced symplectic
structure. Then left translation by $\w^{-1}$ induces a
symplectomorphism from $C_{\w}$, with its image in
$(S(1),\omega)$, which is identified with $\check {N}^{-}\cap
w^{-1}\check {N}^{+}w \subset \check {N}^{-}$.
\end{enumerate}
\end{theorem}

\begin{proof} We claim that we can choose $\w_0$ in Theorem~\ref{bruhatbirkhoff}
so that there is a minimal factorization of the form
$\w_0=r_M..r_{n+1}r_n..r_1$, where each $ r_j$ corresponds to a
simple positive root $\gamma_j$. To prove this it suffices to show
(in the Weyl group) that $N(w_0r_1..r_n)=M-n$, where $N(\cdot)$
denotes the length of a Weyl group element. It suffices to show
that
\begin{equation}\nonumber 
N(w_0r_1..r_s)=N(w_0r_1..r_{s-1})-1
\end{equation}
 for $s=1,..,n$. This is the case precisely when
 $w_0r_1..r_{s-1}\cdot \gamma_s<0$ (see Lemma~4.15.6 of
 \cite{Var}), or equivalently $w_{s-1}^{-1}\cdot \gamma_s>0$. But these
 are precisely the positive roots which are mapped to negative
 roots by $w$ (see Theorem~4.15.10 of \cite{Var}). This proves the claim.

Theorem~\ref{bruhatbirkhoff} asserts that translation by $\w_0$
(or $\w_0^{-1}$) induces a symplectomorphism from the top Bruhat
leaf, with (the negative of) the symplectic structure induced by
$\pi_K$, with $(S(1),\omega )$. We can now directly translate the
results in \cite{Lu} into our framework. When we compose the
parameterization in Theorem~2.1 of \cite{Lu} for the top Bruhat
leaf, with translation by $\w_0^{-1}$, we obtain a
parameterization
$$\C^M\to S(1):(\zeta_M,..,\zeta_1)\to w_{M-1}^{-1}i_{\gamma_
M}(k(\zeta_M))w_{M-1}..w_1^{-1}i_{\gamma_2}(k(\zeta_2))w_1i_{\gamma_
1}(k(\zeta_1)),$$ such that $\omega$, $a(k(\zeta ))$, and
$d\lambda_{N^{-}}$ are expressed as in the statement of the
theorem, with $M$ in place of $n$.

The various parts of the theorem follow from the product structure
of these formulas.
\end{proof}


\begin{thebibliography}{99}
\bibitem[C]{Caine} Caine, A.: \emph{Compact symmetric spaces, triangular factorization, and Poisson geometry}.
arXiv:math/0608454, to appear in J. Lie Th.

\bibitem[EL]{EL} Evens, S; Lu, J.H.: \emph{On the variety of Lagrangian subalgebras, I}.  Ann.
Scient. \'Ec. Norm. Sup. $4^e$ s\'erie, t. 34 (2001) 631-668.

\bibitem[FO]{FO} Foth, P.; Otto, M.: \emph{A symplectic realization of Van Den Ban's convexity
theorem}.  Doc. Math.  11  (2006) 407-424 (electronic).

\bibitem[FL]{FL} Foth, P.; Lu, J.H.: \emph{A Poisson structure on compact symmetric spaces}. Comm. Math. Phys. 251, no. 3 (2004)
557-566.

\bibitem[He]{Helgason} Helgason, S.: Differential Geometry, Lie Groups, and Symmetric Spaces. Pure and Applied Mathematics 80,  Academic Press, New
York-London, 1978.

\bibitem[Lu]{Lu} Lu, J.H.: \emph{Coordinates on Schubert cells,
Kostant's harmonic forms, and the Bruhat-Poisson structure on
$G/B$}. Transf.  Groups 4, No.  4 (1999) 355-374.

\bibitem[Pi1]{pickrell1} Pickrell, D.: \emph{The diagonal distribution for the invariant measure of a
unitary type symmetric space}. Transformation Groups, 11 (4)
(2006) 705-724.

\bibitem[Pi2]{pickrell2} Pickrell, D.: \emph{Homogeneous Poisson
structures on loop spaces of symmetric spaces}, arXiv:0801.3277.

\bibitem[V]{Vaisman} Vaisman, I.: Lectures on the geometry of Poisson manifolds. Progress in Mathematics, 118. BirkhÃ¤user Verlag, Basel,
1994.

\bibitem[Var]{Var} Varadarajan, V.S.:  Lie Groups, Lie Algebras,
 and Their Representations. Springer-Verlag, 1984.

\end{thebibliography}
\end{document}